\documentclass[10pt]{article}
\usepackage{amssymb,amsmath,latexsym,bbm,amscd}

\newtheorem{theorem}[equation]{Theorem}
\newtheorem{corollary}[equation]{Corollary}
\newtheorem{lemma}[equation]{Lemma}
\newtheorem{proposition}[equation]{Proposition}
\newtheorem{definition}[equation]{Definition}

\numberwithin{equation}{section}

\newcommand{\ad}{\hbox{ad}\,}

\newcommand{\Cx}{{\mathbb C}}
\newcommand{\Tx}{{\mathbb T}}

\newcommand{\Zx}{{\mathbb Z}}
\newcommand{\fg}{\mathfrak{g}}
\newcommand{\fh}{\mathfrak{h}}
\newcommand{\aff}{{\hbox{\tiny aff}}}
\newcommand{\fgdot}{{\fg}}
\newcommand{\fgaff}{{\widehat \fg}}

\newcommand{\ot}{\otimes}
\newcommand{\cK}{\mathcal{K}}
\newcommand{\ol}{\overline}
\newcommand{\Z}{\mathbb{Z}}
\newcommand{\cD}{\mathcal{D}}
\newcommand{\Der}{\hbox{Der}\,}

\newcommand{\ga}{\alpha}

\newcommand{\gs}{\sigma}
\newcommand{\wgs}{\widehat{\sigma}}
\newcommand{\bgs}{\overline{\sigma}}
\newcommand{\proof}{{\bf Proof\ \ }}
\newcommand{\qed}{\hfill $\Box$}
\newcommand{\vac}{{\mathbbm 1}}
\newcommand{\fghat}{\widehat{\fg}}

\newcommand{\gd}{\delta}

\newcommand{\cU}{{\mathcal U}}

\newcommand{\cR}{{\mathcal R}}
\newcommand{\cS}{{\mathcal S}}
\newcommand{\cL}{{\mathcal L}}
\newcommand{\cM}{{\mathcal M}}
\newcommand{\cN}{{\mathcal N}}
\newcommand{\cV}{{\mathcal V}}

\newcommand{\FF}{{\mathfrak{F}}}

\newcommand{\td}{{\widetilde d}}

\newcommand{\ou}{{\overline u}}

\newcommand{\oc}{{\overline c}}
\newcommand{\ok}{{\overline k}}
\newcommand{\oal}{{\overline \alpha}}
\newcommand{\obe}{{\overline \beta}}
\newcommand{\orr}{{\overline r}}

\newcommand{\gV}{{\mathfrak {glVir}}}
\newcommand{\sV}{{\mathfrak {slVir}}}
\newcommand{\VglV}{V_{\mathfrak {glVir}}}
\newcommand{\LglV}{L_{\mathfrak {glVir}}}
\newcommand{\omglV}{\omega_{\mathfrak {glVir}}}

\newcommand{\glnhat}{\widehat{{\mathfrak {gl}}}_N}
\newcommand{\gln}{{\mathfrak {gl}}_N}
\newcommand{\sln}{{{\mathfrak {sl}}_N}}
\newcommand{\sll}{{\mathfrak {sl}}}
\newcommand{\so}{{\mathfrak {so}}}

\newcommand{\h}{{\mathfrak h}}
\newcommand{\Heis}{{\tiny {Heis}}}
\newcommand{\hei}{{\tiny {Heis}}}
\newcommand{\vh}{{\tiny {VH}}}
\newcommand{\vir}{{\tiny {Vir}}}
\newcommand{\fgtor}{{\fg_{\hbox{\tiny T}}}}
\newcommand{\VT}{{\cV_{\hbox{\tiny T}}}}
\newcommand{\fgE}{\fg_{\hbox{\tiny E}}}

\newcommand{\End}{\hbox{\rm End\,}}
\newcommand{\Id}{\hbox{\rm Id\,}}

\newcommand{\ols}{\ol{s}}
\newcommand{\p}{\partial}
\newcommand{\dd}{\hbox{d}}

\newcommand{\W}{{\mathcal W}}
\newcommand{\U}{{\mathcal U}}
\newcommand{\Hyp}{V^+_{\hbox{\tiny{Hyp}}}}
\newcommand{\VH}{\Hyp}
\newcommand{\hyp}{{\hbox{\tiny Hyp}}}

\newcommand{\Vaff}{V_{\hbox{\tiny{aff}}}}
\newcommand{\qzz}{\left(\frac{z_2}{z_1}\right)}
\newcommand{\I}{I}

\title{Irreducible Modules for \\ Extended Affine Lie Algebras}

\author{Yuly Billig$\hbox{\,}^{a\,}\hbox{}^*$\ \, and Michael Lau$\hbox{\,}^{b\,}$\thanks{Both authors gratefully acknowledge funding from the Natural Sciences and Engineering Research Council of Canada.}\ \,\footnote{Corresponding author.} \vspace{0.3cm}\\$\hbox{\ \,}^a${\small Carleton University, School of Mathematics and Statistics},\\ {\small 1125 Colonel By Drive, Ottawa, Ontario, Canada K1S 5B6}\\ {\small Email:\ billig@math.carleton.ca}\vspace{0.1cm}\\ $\hbox{\ \,}^b${\small Universit\'e Laval, D\'epartement de math\'ematiques et de statistique},\\ {\small 1045, av. de la M\'edecine, Qu\'ebec,
Qu\'ebec, Canada G1V 0A6}\\ {\small Email:\ Michael.Lau@mat.ulaval.ca}}
\date{}

\begin{document}
\maketitle

\begin{small}
\noindent
{\bf Abstract.}  
We construct irreducible modules for twisted toroidal Lie algebras and extended affine Lie algebras. This is done by combining the representation theory of untwisted toroidal algebras with the technique of thin coverings of modules.  We illustrate our method with examples of extended affine Lie algebras of Clifford type.

\bigskip

\noindent
{\bf MSC: } 17B67, 17B69.

\noindent
{\bf Keywords: }
twisted toroidal Lie algebras, extended affine Lie algebras, thin coverings, twisted modules.

\end{small}

\setcounter{section}{-1}
\vskip.25truein
\section{Introduction}

Extended affine Lie algebras (EALAs) are natural generalizations of the affine Kac-Moody algebras.  They come equipped with a nondegenerate symmetric invariant bilinear form, a finite-dimensional Cartan subalgebra, and a discrete root system.  Originally introduced in the contexts of singularity theory and mathematical physics, their structure theory has been extensively studied for over 15 years.  (See \cite{AABGP,ABFP2,neher} and the references therein.)

Their representations are much less well understood.  Early attempts to replicate the highest weight theory of the affine setting were stymied by the lack of a triangular decomposition; later work considered only the untwisted toroidal Lie algebras and a few other isolated examples.

As a result of major breakthroughs announced in \cite{ABFP2} and \cite{neher}, it is now clear that, except for (extensions of) matrix algebras over non-cyclotomic quantum tori, every extended affine Lie algebra can be constructed as an extension of a twisted multiloop algebra.  These results have inspired the present paper, in which we use a twisting procedure to explicitly obtain irreducible generalized highest weight modules for EALAs associated with every twisted multiloop algebra.

In more detail, let $\fg$ be a finite-dimensional simple Lie algebra over the complex numbers $\Cx$, with commuting automorphisms $\gs_0,\gs_1,\ldots, \gs_N$ of orders $m_0,m_1,\ldots, m_N$, respectively.  Fix primitive $m_i$th roots of unity $\xi_i\in\Cx$ for every $i$, and let each $\gs_i$ act as an automorphism of the $(N+1)$-torus $\Tx^{N+1}=\left(\Cx^{\times}\right)^{N+1}$, by sending the point $(x_0,x_1,\ldots,x_N)\in\Tx^{N+1}$ to the point $(x_0,\ldots,\xi_ix_i\ldots,x_N)$.  The {\em twisted multiloop algebra} $L(\fg;\gs)$ consists of the $\gs_0,\ldots,\gs_N$-equivariant $\fg$-valued regular functions on $\Tx^{N+1}$, under pointwise Lie bracket.  
Next, we take the universal central extension of $L(\fg;\gs)$ and adjoin the Lie algebra of equivariant vector fields on the torus, possibly twisted with a $2$-cocycle.  This produces the {\em full twisted toroidal Lie algebra}
$$\fgtor(\gs) =L(\fg;\gs)\oplus\cK_\Lambda\oplus\cD_\Lambda.$$
This Lie algebra does not admit a nondegenerate invariant symmetric bilinear form, and is thus too large to be an extended affine Lie algebra.  Indeed, the largest extended affine Lie algebra $\fgE(\gs)$ associated with $L(\fg;\gs)$ is the Lie algebra obtained by adjoining only the {\em divergence zero} vector fields to the universal central extension of $L(\fg;\gs)$:  
$$\fgE(\gs) =L(\fg;\gs)\oplus\cK_\Lambda\oplus\cS_\Lambda.$$
See Section 1 for details.


After describing the twisted toroidal Lie algebras, we discuss {\em $\gs$-twists of vertex Lie algebras} in Section 2.  These structures may be thought of as the twisted analogues of the vertex Lie algebras of Dong, Li, and Mason \cite{DLM}.  They are examples of {\em $\Gamma$-twisted formal distribution algebras}, a more general construction appearing in the work of Kac \cite{FDA}.  We go on to prove that the twisted toroidal Lie algebra $\fgtor(\gs_0,1,\ldots,1)$ is a $\gs_0$-twist of the (untwisted) full toroidal Lie algebra $\fgtor(1,\ldots,1)$.  

Past work by one of the authors identifies a quotient $\VT$ of the universal enveloping vertex operator algebra (VOA) of $\fgtor(1,\ldots,1)$ as a tensor product of a lattice VOA, an affine VOA, and a VOA $V_\gV$ associated with the affine Lie algebra ${\widehat{\mathfrak{gl}}_N}$ and the Virasoro Lie algebra.  This allows irreducible representations of $\fgtor(1,\ldots,1)$ and $\fgE(1,\ldots,1)$ to be constructed from tensor products of modules for the tensor components of $\VT$.  See \cite{B1,B2} for details.

In the present paper, we show that a tensor product $\VH\ot\W\ot L_{\gV}$ of modules for the corresponding lattice VOA, twisted affine Lie algebra $\fghat(\gs_0)$, and $V_\gV$ can be given the structure of a module for the Lie algebra $\fgtor(\gs_0,1,\ldots,1)$.  This is done by using general theorems about $\gs_0$-twists of vertex Lie algebras to reduce most of the verifications to work previously done in the untwisted case.

To obtain irreducible modules for the Lie algebras $\fgtor(\gs_0,\gs_1,\ldots,\gs_N)$ and $\fgE(\gs_0,\gs_1,\ldots,\gs_N)$, we use the technique of {\em thin coverings} introduced in a previous paper \cite{BL}. Thin coverings are a tool for constructing graded-simple modules from simple ungraded modules over a graded algebra. By taking a thin covering with respect to $\langle \gs_1\rangle\times\cdots\times\langle \gs_N\rangle$ of an irreducible highest weight module for a twisted affine algebra $\fghat(\gs_0)$, we produce irreducible representations for the twisted toroidal Lie algebras $\fgtor(\gs)$ and $\fgE(\gs)$.  These lowest energy modules have weight decompositions into finite-dimensional weight spaces, and the action of the centres of $\fgtor(\gs)$ and $\fgE(\gs)$ is given by a central character.

We illustrate our method by explicitly constructing irreducible representations for two of the more exotic extended affine Lie algebras.  In the process, we give a detailed discussion of how to realize Jordan torus EALAs of Clifford type as extensions of twisted multiloop algebras, and how to find the thin coverings used in our construction. Vertex operator representations of some Clifford type EALAs were previously constructed in \cite{Ta} and \cite{MT}. Unlike this earlier work, our construction yields irreducible modules.

%

\section{Twisted Toroidal Lie Algebras}
Let $\fgdot$ be a finite-dimensional simple Lie algebra over the complex numbers $\Cx$, with commuting automorphisms $\gs_0,\gs_1,\ldots , \gs_N$ of (finite) orders $m_0,m_1,\ldots ,m_N$, respectively.  Fix primitive $m_i$th roots of unity $\xi_i\in\Cx$ for $i=0,1,\ldots,N$.  Define two sublattices $\Gamma\subseteq\Z^N$ and $\Lambda\subseteq\Z^{N+1}$:
$$\Gamma =  m_1\Z\times\cdots\times m_N\Z, \;\;\;
\Lambda=
m_0\Z\times \Gamma.$$
For each ${\bf s}\in\Z^{N+1}$, we write ${\bf s} = (s_0, s)$, where $s_0\in\Z$, $s\in\Z^N$, and 
denote by $\ol{\bf s} = (\ol{s}_0, \ol{s})$ its image under  the canonical map 
$\Z^{N+1}\rightarrow \Z^{N+1}/\Lambda$.  
Likewise, $f({\bf t})$ will denote a Laurent polynomial in the $N+1$ variables 
$t_0^{\pm 1/{m_0}} ,t_1^{\pm 1},\ldots,t_N^{\pm 1}$.  
However, in recognition of the special role played by the first variable $t_0$, 
the multi-index exponential notation $t^r=t_1^{r_1}t_2^{r_2}\cdots t_N^{r_N}$ 
will be reserved for $N$-tuples $r=(r_1,r_2,\ldots,r_N)\in\Z^N$.

The Lie algebra $\fgdot$ has a common eigenspace decomposition
$$\fgdot=\bigoplus_{\ol{\bf s}\in\Z^{N+1}/\Lambda}\fgdot_{\ol{\bf s}},$$
where $\fgdot_{\ol{\bf s}}=\{x\in\fgdot\ |\ \gs_ix=\xi_i^{s_i}x\hbox{\ for\ }i=0,1,\ldots,N\}$.  The corresponding {\em twisted multiloop algebra}
\begin{align}\label{multiloop}
L(\fg;\gs)&=\sum_{{\bf s}\in\Z^{N+1}}t_0^{{s_0}/{m_0}}t^{s}\ot\fgdot_{\ol{\bf s}}\\
&\subseteq\Cx [t_0^{\pm 1/m_0}, t_1^{\pm 1}, \ldots ,t_N^{\pm 1}]\ot \fgdot
\end{align}
has Lie bracket given by
\begin{equation}\label{multiloopmult}
[f_1({\bf t})g_1,f_2({\bf t})g_2]=f_1({\bf t})f_2({\bf t})[g_1,g_2].
\end{equation}
For simplicity of notation, we sometimes drop the tensor product symbol $\ot$, as in (\ref{multiloopmult}).  

Let $\cR=\Cx[t_0^{\pm 1/m_0}, t_1^{\pm 1}, \ldots ,t_N^{\pm 1}]$ be the algebra of Laurent polynomials,
and let $\cR_\Lambda=\Cx[t_0^{\pm 1}, t_1^{\pm m_1}, \ldots ,t_N^{\pm m_N}]\subseteq \cR$.

We will write $\Omega_{\cR}^1$ (respectively, $\Omega_{\cR_\Lambda}^1$) for the space of K\"ahler differentials of $\cR$ (resp., $\cR_\Lambda$).  
As a left $\cR$-module, $\Omega_{\cR}^1$ has a natural basis consisting of the $1$-forms 
$k_p=t_p^{-1}dt_p$ for $p=0,\ldots,N$.  Likewise,
$$\Omega_{\cR_\Lambda}^1=\bigoplus_{p=0}^N\cR_\Lambda k_p.$$
For each $f\in\cR$, the differential map $\dd:\ \cR\rightarrow\Omega_{\cR}^1$ is defined as
$$\dd(f)=\sum_{p=0}^Nd_p(f)k_p,$$
where $d_p=t_p\frac{\partial}{\partial t_p}$ for $p=0,\ldots ,N$.

Kassel \cite{kassel} has shown that the centre $\cK$ of the universal central extension $(\cR\ot\fgdot)\oplus\cK$ of $\cR\ot\fgdot$ can be realized as
$$\cK=\Omega_\cR^1/\dd(\cR).$$
The multiplication in the universal central extension is given by
\begin{equation}\label{ucemult}
[f_1({\bf t})x,f_2({\bf t})y]=f_1({\bf t})f_2({\bf t})[x,y]+(x\,|\,y)f_2\dd(f_1),
\end{equation}
for all $f_1,f_2\in\cR$ and $x,y\in\fgdot$, where $(x\,|\,y)$ is a symmetric invariant bilinear form on $\fgdot$. This form is normalized by the condition that the induced form on the dual of the Cartan subalgebra satisfies $(\alpha\,|\,\alpha) = 2$ for long roots 
$\alpha$.   
Similarly,  the Lie algebra $L(\fg;\gs)$ can be centrally extended by $\cK_\Lambda=\Omega_{\cR_\Lambda}^1/\dd(\cR_\Lambda)$ using the Lie bracket (\ref{ucemult}). This central extension of the twisted multiloop algebra is also universal \cite{Ne2}.

Let $\cD=\hbox{Der}\,\cR$ be the Lie algebra of derivations of $\cR$, and let $\cD_\Lambda=\hbox{Der}\,\cR_\Lambda\subseteq\cD$.  The space $\cD$ (resp., $\cD_\Lambda$) acts on $\cR\ot\fgdot$ (resp., $L(\fg;\gs)$) by
\begin{equation}
[f_1({\bf t})d_a,f_2({\bf t})x]=f_1 d_a(f_2)x.
\end{equation}
There is also a compatible action of $\cD$ (resp., $\cD_\Lambda$) on $\cK$ (resp., on $\cK_\Lambda$) via the Lie derivative:
\begin{equation}
[f_1({\bf t})d_a,f_2({\bf t})k_b]=f_1 d_a(f_2)k_b + \gd_{ab}f_2 \dd(f_1).
\end{equation}
The multiplication in the semidirect product Lie algebra $(\cR\ot\fgdot)\oplus\cK\oplus\cD$ can be twisted by any $\cK$-valued $2$-cocycle $\tau\in\hbox{H}^2(\cD,\cK)$:
\begin{equation}\label{cocyclemult}
[f_1({\bf t})d_a,f_2({\bf t})d_b]=f_1d_a(f_2)d_b-f_2d_b(f_1)d_a+\tau(f_1d_a,f_2d_b).
\end{equation}
We will use cocycles
\begin{equation}\label{cocycleparams}
\tau=\mu\tau_1+\nu\tau_2
\end{equation}
parametrized by $\mu,\nu\in\Cx$.  To define these cocycles, recall that the {\em Jacobian} $v^J$ of a vector field $v=\sum_a f_a({\bf t}) d_a$ is the matrix with $(a,b)$-entry $d_b(f_a)$, for $0\leq a,b\leq n$.  In this notation,
\begin{eqnarray*}
\tau_1(v,w)&=&Tr(v^J\dd(w^J))\\
\tau_2(v,w)&=&Tr(v^J)\dd(Tr(w^J)),
\end{eqnarray*}
where $Tr$ denotes the trace and the differential map $\dd$ is defined element-wise on the matrix $w^J$.

The resulting Lie algebra
\begin{equation}
\fgtor=(\cR\ot\fgdot)\oplus\cK\oplus_\tau\cD
\end{equation}
is called the {\em  toroidal Lie algebra}.  When restricted to $\cD_\Lambda$, the cocycle $\tau$ restricts to a cocycle (also denoted by $\tau$) in the space $\hbox{H}^2(\cD_\Lambda,\cK_\Lambda)$.  This gives the {\em full twisted toroidal Lie algebra}
\begin{equation}
\fgtor(\gs)=L(\fg;\gs)\oplus\cK_\Lambda\oplus_\tau\cD_\Lambda
\end{equation}
with Lie bracket given by (\ref{ucemult}) -- (\ref{cocyclemult}).

We will also consider the closely related (twisted) toroidal extended affine Lie algebra 
(EALA).  A derivation $v$ is called {\em divergence zero} 
(or {\em skew-centroidal}) if $Tr(v^J)=0$.  We will denote the subalgebra of divergence zero derivations of $\cR$ (resp., $\cR_\Lambda$) by $\cS$ (resp., $\cS_\Lambda$). Note that the cocycle $\tau_2$ vanishes when restricted to the space $\cS$, so when working in the EALA setting, we can assume that $\tau=\mu\tau_1$.  The {\em toroidal EALA} is the Lie algebra
\begin{equation}
\fgE=(\cR\ot\fgdot)\oplus\cK\oplus_\tau\cS \subseteq\fgtor.
\end{equation}
Analogously, the {\em twisted toroidal EALA} is the subalgebra
\begin{equation}
\fgE(\gs)=L(\fg;\gs)\oplus\cK_\Lambda\oplus_\tau\cS_\Lambda \subseteq \fgtor(\gs).
\end{equation}
The Lie algebras $\fgE$ and $\fgE(\gs)$ possess non-degenerate invariant bilinear forms.

\section{Vertex Lie Algebras and their $\gs$-Twists}

In this section, we describe a general construction that will reduce the work of verifying certain relations in the twisted toroidal setting to verifying the analogous relations in the untwisted setting.  We begin by recalling the definition of {\em vertex Lie algebra}.  In our exposition, we will follow the paper of Dong, Li, and Mason \cite{DLM}.  Similar constructions appear in the work of Kac \cite{VA}, under the name {\em Lie formal distribution algebra}.  

Let $\mathcal{L}$ be a Lie algebra with basis $\{u(n), c(-1)\ |\ u\in\mathcal{U},c\in \mathcal{C}, n\in\Z\}$, where $\mathcal{U}$ and $\mathcal{C}$ are some index sets. Define the corresponding formal fields in $\mathcal{L}[[z,z^{-1}]]$:
\begin{align}
u(z)&=\sum_{n\in\Z} u(n)z^{-n-1}\\
c(z)&=c(-1)z^0,
\end{align}
for each $u\in\cU$ and $c\in\mathcal{C}$.
Let $\mathcal F$ be the subspace of $\cL[[z,z^{-1}]]$ spanned by the fields $u(z),c(z)$, and their derivatives of all orders.

The {\em delta function} is defined as
$$\gd(z)=\sum_{j\in\Z}z^j .$$

\begin{definition}{\em A Lie algebra $\cL$ with basis as above is called a {\em vertex Lie algebra} if the following two conditions hold:
\begin{enumerate}
\item[{\rm (VL1)}] For all $u_1,u_2\in\cU$, there exist $n\geq 0$ and $f_0(z),\ldots,f_n(z)\in\mathcal{F}$ such that
$$[u_1(z_1),u_2(z_2)]=\sum_{j=0}^nf_j(z_2)\left[ z_1^{-1}\left(\frac{\partial}{\partial z_2}\right)^j\gd\left(\frac{z_2}{z_1}\right)\right].$$
\item[{\rm (VL2)}] The element $c(-1)$ is central in $\cL$ for all $c\in\mathcal{C}$.
\end{enumerate}
}
\end{definition}
Let $\cL^{(-)}$ be the subspace with basis $\{u(n),c(-1)\ |\ u\in\mathcal{U},\ c\in\mathcal{C},\ n<0\}$, and let $\cL^{(+)}$ be the subspace of $\cL$ with basis $\{u(n)\ \ u\in\mathcal{U},\ n\geq 0\}$.  Then $\cL=\cL^{(-)}\oplus\cL^{(+)}$ and both $\cL^{(-)}$ and $\cL^{(+)}$ are in fact subalgebras of $\cL$.

The universal enveloping vertex algebra $V_\cL$ of a vertex Lie algebra $\cL$ is the induced module
\begin{equation}
V_\cL=\hbox{Ind}_{\cL^{(+)}}^{\cL}(\Cx \vac)=U(\cL^{(-)})\ot_\Cx\vac,
\end{equation}
where $\Cx\vac$ is a trivial $1$-dimensional $\cL^{(+)}$-module.

The following result appears as Theorem 4.8 in \cite{DLM}.  (See also \cite{VA}.)
\begin{theorem}\label{reconstruction} Let $\cL$ be a vertex Lie algebra.  
Then $V_\cL$ has the structure of a vertex algebra with vacuum vector $\vac$. 
The infinitesimal translation operator $T$ is the derivation of 
$V_\cL$ given by $T(u(n))=$ 
\break
$-nu(n-1)$ and $T(c(-1))=0$ for all $u\in\mathcal{U}$, $c\in\mathcal{C}$.  The state-field correspondence is defined by the formula
\begin{align*}
&Y\big(a_1(-n_1-1)\cdots a_{k-1}(-n_{k-1}-1)a_k(-n_k-1)\vac,z\big)\\
&\ \ =:\left(\frac{1}{n_1!}\left(\frac{\partial}{\partial z}\right)^{n_1}a_1(z)\right)\cdots\\
&\ \ \ \ \ \ :\left(\frac{1}{n_{k-1}!}\left(\frac{\p}{\p z}\right)^{n_{k-1}}a_{k-1}(z)\right)\left(\frac{1}{n_k!}\left(\frac{\p}{\p z}\right)^{n_k}a_k(z)\right):\cdots :,
\end{align*}
where $a_j\in\mathcal{U}$ and $n_j\geq 0$, or $a_j\in\mathcal{C}$ and $n_j=0$.
\end{theorem}

Next we will define a twisted vertex Lie algebra (cf. \cite{FDA,KW}).  

\bigskip

We consider a vertex Lie algebra $\cL$ graded by a cyclic group $\Z/m\Z$ for which the generating fields $u(z)$, $c(z)$ are homogeneous.  
That is, $\cL=\bigoplus_{\ol{k}\in\Z/m\Z}\cL_{\ol{k}}$ is a $\Z/m\Z$-graded Lie algebra, and there is a decomposition 
$\mathcal{U}=\bigcup_{\overline{k}\in\Z/m\Z}\mathcal{U}_{\overline{k}}$ for which $u(n)\in\cL_{\overline{k}}$  and $c(-1)\in\cL_{\overline{0}}$  
for all $c\in\mathcal{C}$, $u\in\mathcal{U}_{\overline{k}}$, and $n\in\Z$.  
Let $\xi$ be a primitive $m$th root of 1.
This grading defines an automorphism $\gs:\ \cL\rightarrow\cL$ 
of order $m$ by $\gs(x) = \xi^k x$ for $x \in \cL_{\ol{k}}$.

Let $\cL(\gs)$ be a space with the basis
 \begin{equation}\label{Lsigmadef}
\bigcup_{\ok\in\Z/m\Z}\left\{\ou(n),\oc(-1)\ \left|\ u\in\mathcal{U}_{\overline{k}},\ n\in {k}/{m}+\Z,\ c\in\mathcal{C}\right.\right\}.
\end{equation}

We define fields
\begin{eqnarray}
\overline{u}(z)&=&\sum_{n\in {k/m}+\Z}\ou(n)z^{-n-1}\in\cL(\gs)[[z^{-1/m},z^{1/m}]],\\
\overline{c}(z)&=&\oc(-1)z^0,
\end{eqnarray}
for all $u\in\mathcal{U}_{\overline{k}}$ and $c\in\mathcal{C}$.  Let $\overline{\mathcal{F}}$ be the space spanned by the fields $\overline{u}(z),\overline{c}(z)$, and their derivatives of all orders.  
The correspondence $u(z)\mapsto\overline{u}(z)$, $c(z)\mapsto\overline{c}(z)$ extends to a vector space isomorphism 
$\overline{\hbox{\; \vbox{\vspace{0.2cm}}}}: \;\; \mathcal{F}\rightarrow\overline{\mathcal{F}}$ commuting with the derivative $\frac{d}{dz}$. 

We will use the {\em twisted delta function} $\delta_k(z)=z^{k/m} \delta (z)$ when working with $\gs$-{\em twists} of vertex Lie algebras.  More precisely, the {\em $\gs$-twist of a vertex Lie algebra} $\cL$ is a vector space $\cL(\gs)$ equipped with a Lie bracket defined by the relations
\begin{enumerate}
\item[{\rm (T1)}]$\displaystyle{[\ol{u}_1(z_1),\ol{u}_2(z_2)]=\sum_{j=0}^n\ol{f}_j(z_2)\left[z_1^{-1}\left(\frac{\p}{\p z_2}\right)^j\gd_k\left(\frac{z_2}{z_1}\right)\right]},$
\item[]where the $f_j$ are as in (VL1) and $u_1\in\mathcal{U}_{\ol{k}}$, $u_2\in\mathcal{U}$.
\item[{\rm (T2)}] The elements $\ol{c}(-1)$ are central in $\cL(\gs)$ for all $c\in\mathcal{C}$,
\end{enumerate}
It follows from \cite{KW} that $\cL(\gs)$ is indeed a Lie algebra.

Observe that the twisted affine Kac-Moody algebras are examples of twisted vertex Lie algebras.  In this paper, our main example of a vertex Lie algebra will be the full toroidal Lie algebra $\fgtor$.  Its $\gs_0$-twist will be the twisted toroidal Lie algebra $\fgtor (\gs_0, 1, \ldots, 1)$.

The next theorem (see e.g., \cite{li}), will be very helpful in our construction of modules for twisted toroidal Lie algebras.  

\begin{theorem}\label{twistedmodulethm}
Let $\cL(\gs)$ be a $\gs$-twist of a vertex Lie algebra $\cL$, and let $V_\cL$ be the universal enveloping vertex algebra of $\cL$.  
Then every $\gs$-twisted $V_\cL$-module $M$ is a module for the Lie algebra $\cL(\gs)$.
\end{theorem}

 Let us recall the definition of a twisted module of a vertex algebra \cite{li,KR}. Let $\gs$ be an automorphism of order $m$ of a vertex algebra $V$. Consider the grading of $V$ by the cyclic group $(\frac{1}{m}\Z)/\Z$, where for each coset $\oal = k/m +\Z$ we define the component
$$ V_\oal = \left\{ v \in V\ |\ \gs(v) = \xi^k v \right\} .$$

For each coset $\oal = k/m +\Z$, fix a representative $\alpha \in \oal$. If $\oal = \Z$, we set $\alpha = 0$.  Let $M$ be a vector space with a map
\begin{equation}
Y_M: \;\; V \rightarrow \End(M)[[z^{1/m},z^{-1/m}]] .
\end{equation}
 Write $Y_M(a,z) = \sum_{j\in \oal} a^M_{(j)} z^{-j-1}$, with each $a^M_{(j)} \in \End (M)$.

\begin{definition}
A vector space $M$ together with a map $Y_M$ as above
is called a $\gs$-twisted module for $V$ if the following axioms hold:
\begin{equation}
Y_M(a,z) \in z^{-\alpha} \End(M)[[z,z^{-1}]], 
\end{equation}
\begin{equation}
a^M_{(\alpha+n)} v = 0 \hbox{\ for \ } n\gg 0,
\end{equation}
\begin{equation}
Y_M(\vac,z) = \Id_M z^0,
\end{equation}
\begin{equation}
Y_M(T(a), z) = \frac{d}{dz} Y_M(a,z),
\end{equation}
\begin{eqnarray}
&{\hskip 0.8cm} 
\sum\limits_{j=0}^\infty {{m}\choose{j}} Y_M (a_{(n+j)} b,z ) z^{m-j} {\hskip 0.6cm} 
\nonumber\\
&= \sum\limits_{j=0}^\infty (-1)^j  {{n}\choose{j}}\label{twistedBorcherds}
\left( a^M_{(m+n-j)} Y_M(b,z) z^j 
- (-1)^n Y_M(b,z) a^M_{(m+j)} z^{n-j} \right) {\hskip 1cm}
\end{eqnarray}
$\hbox{\ for all \;} a\in V_\oal, \,  m \in \alpha+\Z, \, b \in V, v\in M,\hbox{\ and\ } n\in\Z .$

\end{definition}

Letting $n=0$ in the twisted Borcherds identity (\ref{twistedBorcherds}), one gets the commutator formula
for $a \in V_\oal, b \in V_\obe, m\in\oal,$ and $n\in\obe$\,:
\begin{equation}
[ a^M_{(m)}, b^M_{(n)} ] = \sum\limits_{j=0}^\infty {{m}\choose{j}} \left( a_{(j)} b\right)^M_{(m+n-j)}.
\end{equation}

The twisted normally ordered product is defined as (see \cite{KW})
\begin{equation}
:Y_M(a,z) Y_M(b,z) : \; = Y_M(a,z)_+ Y_M(b,z) + Y_M(b,z) Y_M(a,z)_- \; ,
\end{equation}
where
\begin{equation}
Y_M(a,z)_- = \sum\limits_{j=0}^\infty a^M_{(\alpha+j)} z^{-\alpha-j-1}, \;
Y_M(a,z)_+ = \sum\limits_{j=-1}^{-\infty} a^M_{(\alpha+j)} z^{-\alpha-j-1}.
\end{equation}
Note that when $a \in V_{\overline 0}$, this coincides with the usual normal ordered product.

Letting $m=\alpha$ and $n=-1$ in the twisted Borcherds identity, one gets:
\begin{equation}{\label {TBI}}
\sum\limits_{j=0}^\infty {{\alpha}\choose{j}} Y_M (a_{(-1+j)} b, z ) z^{-j} = :Y_M (a,z) Y_M(b,z): .
\end{equation}

\section{Representations of Twisted Toroidal Lie Algebras}

 In this section, we use the representation theory of the full toroidal Lie algebras to construct
irreducible representations for the full twisted toroidal Lie algebras.  We begin by describing the toroidal vertex operator algebra (VOA) that controls the representation theory
of the full toroidal Lie algebra $\fgtor$. We then show that twisted modules for this toroidal VOA
yield representations of the twisted toroidal Lie algebra $\fgtor(\gs_0, 1, \ldots, 1)$. Finally, we realize irreducible modules for $\fgtor(\gs_0, \gs_1, \ldots, \gs_N)$ as subspaces of the
irreducible modules for $\fgtor(\gs_0, 1, \ldots, 1)$ using thin coverings.

\subsection{The toroidal VOA $\VT$}
 The toroidal vertex operator algebra  $\VT$ that controls the representation theory of the full toroidal Lie algebra $\fgtor$ is
a quotient of the universal enveloping vertex algebra of $\fgtor$\ \cite{B1}.  It is a tensor product of three VOAs:
\begin{equation}\label{tortensor}
\VT = \VH \otimes \Vaff \otimes \VglV.
\end{equation}
Here $\VH$ is a sub-VOA of a lattice VOA, $\Vaff$ is an affine VOA, and $\VglV$ is a twisted
$\glnhat$-Virasoro VOA. 
We will give brief descriptions of these VOAs and refer to \cite{B1} for details.

The vertex operator algebra $\VH$ is a sub-VOA of a lattice vertex algebra associated with a hyperbolic lattice.
As a vector space, it is a tensor product of a Laurent polynomial algebra with a Fock
space: 
$$\VH = \Cx[q_1^{\pm 1},\ldots, q_N^{\pm 1}] \otimes \FF,$$
$$\FF = \Cx[u_{pj},v_{pj}\ |\ p=1,2,\ldots,N, \; j = 1, 2, 3, \ldots ].$$
In the description of the action of $\fgtor$ on $\VT$, we will use the following vertex operators:
\begin{eqnarray*}
K_0(r,z) &=& Y(q^r, z)\\
 &=& q^r\exp\left(\sum_{p=1}^Nr_p\sum_{j=1}^{\infty} u_{pj}z^j\right)
\exp\left(-\sum_{p=1}^Nr_p\sum_{j=1}^\infty\frac{z^{-j}}{j}\frac{\p}{\p v_{pj}}\right),\\
K_a(z)&=& Y(u_{a1},z) = \sum_{j=1}^\infty ju_{aj}z^{j-1}+\sum_{j=1}^\infty \frac{\p}{\p v_{aj}}z^{-j-1},\\
K_a(r,z)&=&Y(u_{a1}q^r,z) = K_a(z)K_0(r,z),\\
D_a(z)&=& Y(v_{a1},z) =\sum_{j=1}^{\infty} jv_{aj}z^{j-1}+q_a\frac{\p}{\p q_a}z^{-1}+
\sum_{j=1}^{\infty} \frac{\p}{\p u_{aj}}z^{-j-1},\\
\omega_\hyp (z) &=& Y \left( \sum_{p=1}^N u_{p1} v_{p1}, z\right) = \sum_{p=1}^N :K_p(z) D_p(z): \; ,\\
\end{eqnarray*}
for $a=1,2,\ldots, N$ and $r\in\Z^N$.  The last expression is the Virasoro field of this VOA, and the rank of $\VH$ is $2N$.  

\medskip

\noindent
{\bf Remark.} The vertex algebra $\VH$ has a family of modules
\begin{equation}\label{MHypdef}
M_\hyp (\alpha,\beta) = e^{\beta v}q^\alpha \Cx[q_1^{\pm 1},\ldots, q_N^{\pm 1}] \otimes
\FF,
\end{equation}
where $\ga\in\Cx^N$, $\beta\in\Z^N$, and $\beta v= \beta_1v_1+\cdots +\beta_Nv_N$.  
See \cite{B1} for the description of the action of $\VH$ on $M_\hyp (\alpha,\beta)$. 
All constructions of 
modules in this paper admit a straightforward generalization by shifting the algebra of Laurent polynomials by the
factor $e^{\beta v}q^\alpha$ and replacing $Y(q^r,z)$ with $Y_{M_\hyp(\alpha,\beta)}(q^r,z)$.
For the sake of simplicity of exposition, we will not be using these modules in the present paper.

The second factor in (\ref{tortensor}) is the usual affine vertex operator algebra $\Vaff$ of noncritical level $c\in\Cx$ associated with the affine Lie algebra\begin{equation*}
\fgaff = \big(\fg\ot\Cx[t_0,t_0^{-1}]\big)\oplus\Cx C_\aff.
\end{equation*}
We denote its Virasoro field by $\omega_\aff(z)$.  The corresponding Virasoro algebra has 
central charge $\displaystyle{{c\dim\fg}/{(c+h^\vee)}}$, where $h^\vee$ is the dual Coxeter 
number of $\fg$.  This material may be found in any of the introductory books on vertex 
operator algebras.  (See \cite{VA}, for instance.)

The remaining VOA in the tensor product (\ref{tortensor}) is associated with the
{\em twisted $\glnhat$-Virasoro algebra} $\gV$, which is the universal central extension of
the Lie algebra
\begin{equation}
\big(\Cx[t_0, t_0^{-1}]\ot\gln(\Cx)\big)\rtimes\Der\Cx[t_0, t_0^{-1}].
\end{equation}
This central extension is obtained by adjoining a $4$-dimensional space spanned by the
basis $\{C_{\sln},C_\Heis,C_\vir,C_\vh \}.$

We fix the natural projections
\begin{eqnarray} \label{proj} 
\psi_1:\ \gln(\Cx)&\rightarrow& \sln(\Cx)\\
\psi_2:\ \gln(\Cx)&\rightarrow& \Cx
\end{eqnarray}
where $\psi_2(u)= Tr(u)/N$, $\psi_1(u)=u-\psi_2(u)I$, and $I$ is the $N\times N$ identity
matrix.
The multiplication in $\gV$ is given by
\begin{eqnarray*}
\lbrack L(n),L(m)\rbrack&=&(n-m)L(n+m)+\frac{n^3-n}{12}\gd_{n+m,0}C_\vir\\
\lbrack L(n),u(m)\rbrack&=&-mu(n+m)-(n^2+n)\gd_{n+m,0}\psi_2(u)C_\vh\\
\lbrack u(n),v(m)\rbrack&=&[u,v](n+m)\\
&&\ \ \ \ \ \ \ \ +n\gd_{n+m,0}\big( Tr \,(\psi_1(u)\psi_1(v))C_{\sln}+\psi_2(u)\psi_2(v)C_\Heis\big),
\end{eqnarray*}
where $L(n)$ is the Virasoro operator $-t_0^{n+1}\frac{\p}{\p t_0}$ and $u(m)=t_0^m\ot u$
for $u\in\gln(\Cx)$.

The twisted $\glnhat$-Virasoro algebra $\gV$ is a vertex Lie algebra \cite[Prop 3.5]{B1}, and let $\VglV$ be its universal enveloping vertex algebra with central charge given by a central character $\gamma$.  We write
\begin{eqnarray*}
c_{\sln}&=&\gamma(C_{\sln}),\\
c_\Heis&=&\gamma(C_\Heis),\\
c_\vir&=&\gamma(C_\vir),\\
c_\vh&=&\gamma(C_\vh).
\end{eqnarray*}
The Virasoro field of $\VglV$ is
$$\omglV (z) = Y(L(-2) \vac, z) = \sum_{j\in\Z} L(j) z^{-j-2}.$$

For $i,j=1,\ldots, N$, let 
$$E_{ij} (z) = Y(E_{ij}(-1) \vac , z) = \sum_{k\in\Z} E_{ij}(k) z^{-k-1},$$
where $E_{ij} \in \gln$ is the matrix with $1$ in the $(i,j)$-position, and zero elsewhere.

\subsection{Representations of $\fgtor(\gs_0,1,\ldots,1)$}
For $a=1,\ldots,N$ and ${\bf r} = (r_0, r) \in\Z^{N+1}$, 
we now define fields (in a single variable $z$) whose Fourier coefficients span the Lie algebra $\fgtor(\gs_0,1,\ldots, 1)$:
\begin{eqnarray}\label{field1eqn}
k_0(r,z)&=&\sum_{j\in\Z}t_0^{j}t^rk_0z^{-j},\\
k_a(r,z)&=&\sum_{j\in\Z}t_0^{j}t^rk_az^{-j-1},\\
x(r,z)&=&\sum_{j\in{{r_0}/{m_0}}+\Z}t_0^{j}t^rxz^{-j-1}\hbox{\ for each\ }x\in\fgdot_{\ol{\bf r}},\\
\widetilde{d}_a(r,z)&=&\sum_{j\in\Z}\big(t_0^{j}t^rd_a-\nu r_at_0^{j}t^rk_0\big)z^{-j-1},\\
\widetilde{d}_0(r,z)&=&-\sum_{j\in\Z}\left(t_0^{j}t^rd_0-(\mu+\nu)\left(j+\frac12\right)t_0^{j}t^rk_0\right)z^{-j-2},\label{lastfieldeqn}
\end{eqnarray}
where $\mu$ and $\nu$ are the parameters of the cocycle $\tau=\mu\tau_1+\nu\tau_2$, as in (\ref{cocycleparams}).
In the case where $\gs_0=1$ (and $m_0=1$), these fields may be viewed as the generating fields of the (untwisted) full toroidal Lie algebra $\fgtor$.  (See \cite{B1} for details.)

We consider the commutation relations between these fields.  Most of these relations can be 
taken directly from work done for the untwisted case \cite[Eqn 5.7]{B1}.  The only exceptions are those relations involving fractional powers of $z$--namely, the relations involving the field $x(r,z)$.  Verifying these relations is a completely straightforward calculation.  For $x\in\fg_{\ol{\bf r}}$, $y\in\fg_{\ol{\bf s}}$, $a=1,\ldots, N$, and $i=0,\ldots, N$, we see that

\begin{eqnarray}\label{f}
[x(r,z_1),y(s,z_2)]&=&[x,y](r+s,z_2)\left[z_1^{-1}\gd_{r_0}\qzz\right] \nonumber\\
&&+(x|y)k_0(r+s,z_2)\left[z_1^{-1}\frac{\p}{\p z_2} \gd_{r_0}\qzz \right] \nonumber\\
&&+(x|y)\sum_{p=1}^Nr_pk_p(r+s,z_2)\left[z_1^{-1}\gd_{r_0}\qzz\right],\\
\lbrack\widetilde{d}_a(r,z_1),y(s,z_2)\rbrack&=&s_ay(r+s,z_2)\left[z_1^{-1}\gd\qzz\right],\\
\lbrack\widetilde{d}_0(r,z_1),y(s,z_2)\rbrack&=&\frac{\p}{\p z_2}\left(y(r+s,z_2)\left[z_1^{-1}\gd\qzz\right]\right),\\
\lbrack x(r,z_1),k_i(s,z_2)\rbrack&=&0. \label{l}
\end{eqnarray}

The above computations demonstrate that the twisted toroidal Lie algebra $\fgtor(\gs_0,1,\ldots, 1)$ is a $\gs_0$-twist of  the untwisted toroidal Lie algebra $\fgtor$. Indeed, the automorphism
$\gs_0:\fg\rightarrow\fg$ naturally lifts to an automorphism
\begin{equation}\label{gs0lift}
\gs_0:\ \fgtor \rightarrow\fgtor
\end{equation}
by setting $\gs_0( t_0^{r_0} t^r x )=t_0^{r_0} t^r \gs_0(x) =\xi_0^{r_0}t_0^{r_0} t^r x $  for each $x\in\fg_{\ol{\bf r}}$ and letting $\gs_0$ act trivially on $\cD$ and $\cK$.
 Comparing (\ref{f})--(\ref{l}) with the corresponding computations in the untwisted toroidal case \cite[Eqn. 5.7]{B1}, we have now verified the following proposition:
\begin{proposition}
The Lie algebra $\fgtor(\gs_0,1,\ldots,1)$ is a $\gs_0$-twist of the vertex Lie algebra $\fgtor$.\qed
\end{proposition}

It was shown in \cite{B1} that for certain central characters $\gamma$, the toroidal VOA $\VT=\VH \otimes \Vaff \otimes \VglV$ is a quotient of the universal enveloping vertex algebra $V_{\fgtor}$.
We will see that the automorphism $\gs_0$ induces an automorphism (again denoted by $\gs_0$) of $\VT$ that is compatible with the natural lift of $\gs_0$ to $V_{\fgtor}$.  Then every $\gs_0$-twisted $\VT$-module is a $\gs_0$-twisted $V_{\fgtor}$-module, and also a $\fgtor(\gs_0,1,\ldots,1)$-module by Theorem \ref{twistedmodulethm}.


\begin{theorem} {\em {\cite{B1}}}
(i) Let $\Vaff$ be the universal enveloping vertex algebra for $\fgaff$ at nonzero, non-critical level $c$. Let $\VglV$ be the universal enveloping vertex algebra of the Lie algebra
$\gV$ with the following central character:
\begin{eqnarray}
c_\sln = 1 - \mu c, \quad c_\hei = N(1-\mu c) - N^2 \nu c,  \nonumber\\
c_\vh = N(\frac{1}{2} - \nu c), \quad
c_\vir = 12c (\mu + \nu) - 2N - \frac{c\dim\fg}{c+h^\vee}, \label{chargammaii}
\end{eqnarray}
where $\mu$ and $\nu$ are as in (\ref{cocycleparams}). 
Then there exists a homomorphism of vertex algebras
$$\phi: \; V_\fgtor \rightarrow  \VH \otimes \Vaff \otimes \VglV ,$$
defined by the correspondence of fields:
\begin{eqnarray}\label{image1eqn}
k_0(r,z) &\mapsto& cK_0(r,z), \\
k_a(r,z) &\mapsto& cK_a(r,z), \label{imageiieqn} \\
x(r,z) &\mapsto& Y(x(-1)\vac, z) K_0(r,z), \label{ximageeqn}\\
\td_a (r,z) &\mapsto& :D_a(z) K_0(r,z): + \sum_{p=1}^N r_p E_{pa}(z) K_0(r,z), \label{imageveqn}\\
\td_0 (r,z) &\mapsto& :\bigg(\omega_\hyp(z) + \omega_\aff(z) + \omglV (z) \bigg) K_0(r,z): \nonumber\\
&&+ \sum_{1\leq i,j\leq N} r_i K_j(z) E_{ij}(z) K_0(r,z) \nonumber \\
&&+ (\mu c - 1) \sum_{p=1}^N r_p \left( \frac{\partial}{\partial z} K_p(z) \right) K_0(r,z),\label{lastimageeqn}
\end{eqnarray}
for all $r\in\Z^N$, $x\in\fg$, and $a=1,\ldots,N$.

(ii) Let $\U$ be an irreducible $\Vaff$-module of nonzero, non-critical level $c$ and let 
$\LglV$ be an irreducible $\VglV$-module with the above central character. Then 
(\ref{image1eqn})-(\ref{lastimageeqn}) define the structure of an irreducible $\fgtor$-module on
$$\Cx [q_1^{\pm 1}, \ldots, q_N^{\pm 1}] \ot \FF \ot \U \ot \LglV .$$
\qed
\end{theorem}

\bigskip

The automorphism $\gs_0:\fg\rightarrow\fg$ lifts to an automorphism of $\fgtor$ as in (\ref{gs0lift}).  It then lifts in the obvious way to a VOA automorphism $\gs_0:\ V_{\fgtor}\rightarrow V_{\fgtor}$.  
It can also be extended to an automorphism $\fgaff \rightarrow \fgaff$ and to
$\gs_0:\ V_\aff\rightarrow V_\aff$ by setting
$$\gs_0(f(t_0) x )=f(t_0) \gs_0(x),\quad \gs_0(C_\aff)=C_\aff,$$
for all $x\in\fg$.  This lets us identify $\gs_0$ with the map $1\ot\gs_0\ot 1$ on the tensor product $\VT$:
$$1\ot\gs_0\ot 1:\ \VH \otimes \Vaff \otimes \VglV\rightarrow\VH \otimes \Vaff \otimes \VglV.$$

\begin{lemma}\label{3.24}
The automorphisms $\gs_0$ on $V_{\fgtor}$ and $\VT$ are compatible with the homomorphism $\phi:\ V_{\fgtor}\rightarrow\VT$ in the sense that
$$\gs_0\circ\phi=\phi\circ\gs_0.$$
\end{lemma}

\noindent
\proof Note that the only field of (\ref{field1eqn})--(\ref{lastfieldeqn}) that is not fixed by $\gs_0$ is $x(r,z)$:
$$\gs_0(x(r,z))=\xi_0^{r_0}x(r,z)$$
for all $x\in\fg_{\ol{\bf r}}$.  The action of $\gs_0$ on the right-hand sides of (\ref{image1eqn})--(\ref{lastimageeqn}) is clearly also trivial, with the exception of its action on (\ref{ximageeqn}) and (\ref{lastimageeqn}).  Note that the only term of (\ref{lastimageeqn}) on which $\gs_0$ can act nontrivially is the affine Virasoro field
$$\omega_\aff(z)=\frac{1}{2(c+h^\vee)}\sum_{i=1}^{\dim\fg}:x_i(z)x^i(z):,$$
where $\{x_i\}$ and $\{x^i\}$ are dual bases of $\fg$, relative to the normalized invariant bilinear form.  These bases can be chosen to consist of eigenvectors for $\gs_0$. 
 Since the invariant bilinear form is $\gs_0$-invariant
(Lemma \ref{invariantform} in Appendix below), we see that the product of eigenvalues of $x_i$ and $x^i$ is $1$, for each $i$.  This means that $\omega_\aff(z)$ is also fixed by $\gs_0$.  On the remaining term (\ref{ximageeqn}), it is clear that for homogeneous $x\in\fg$, the eigenvalues of $\gs_0$ agree on the left- and right-hand sides.  Therefore,
$\phi\circ\gs_0=\gs_0\circ\phi$.\qed

\bigskip

Let $\W$ be an irreducible $\gs_0$-twisted module of the affine vertex operator algebra $V_\aff$.  
The following lemma says that the tensor product $\VH  \otimes \W \otimes \LglV$ is a $\gs_0$-twisted module of $\VT=\VH \otimes \Vaff \otimes \VglV$.

\begin{lemma}\label{3.25}
Let $A$ and $B$ be vertex operator algebras, and assume that $A$ is equipped with a finite-order automorphism $\eta$.  Extend $\eta$ to $A\ot B$ as $\eta\ot 1:\ A\ot B\rightarrow A\ot B$.  If $U$ is an $\eta$-twisted module for $A$ and $V$ is a module for $B$, then $U\ot V$ is an $\eta$-twisted module for $A\ot B$.
\end{lemma}

\noindent
\proof This lemma is a straightforward consequence of \cite[Prop. 3.17]{li}.
\qed

\bigskip

\begin{corollary} \label{1twistmodule}  Let $\W$ be an irreducible $\gs_0$-twisted $V_\aff$-module.  Then the tensor product 
$\VH \otimes \W \otimes \LglV$ is a module for the Lie algebra $\fgtor(\gs_0,1,\ldots,1)$.
\end{corollary}

\noindent
\proof By Lemma \ref{3.25}, $\VH \otimes \W \otimes \LglV$ is a $\gs_0$-twisted module for $\VT$.  By Lemma \ref{3.24}, 
every $\gs_0$-twisted $\VT$-module is also a $\gs_0$-twisted $V_{\fgtor}$-module.  Finally, since $\fgtor(\gs_0,1,\ldots,1)$ 
is a $\gs_0$-twist of $\fgtor$, we see that $\VH \otimes \W \otimes \LglV$ is a module for $\fgtor(\gs_0,1,\ldots,1)$, 
by Theorem \ref{twistedmodulethm}.\qed

\bigskip

Note that by Theorem \ref{twistedmodulethm}, $\gs_0$-twisted $V_\aff$-modules are bounded modules for the twisted affine Lie algebra 
$$\fgaff(\gs_0) = \sum_{j\in\Z} t_0^{j/{m_0}} \fg_{\overline j} \oplus \Cx C_\aff .$$

In order to explicitly describe the action of the Lie algebra $\fgtor(\gs_0,1,\ldots,1)$,
we need to modify formulas (\ref{ximageeqn}) and (\ref{lastimageeqn}) in (\ref{image1eqn})--(\ref{lastimageeqn}).
  
The twisted field $x(r,z)=\sum_{j\in{{r_0}/{m_0}}+\Z}t_0^{j}t^rxz^{-j-1}$ for each 
$x\in\fgdot_{\ol{\bf r}}$, is represented by the twisted vertex operator
\begin{equation}
Y_\W (x(-1)\vac, z) K_0(r,z), \label{twistedaff}
\end{equation}
where $Y_\W (x(-1)\vac, z)$ represents the action of the twisted affine field
$x(z)=\sum_{j\in{{r_0}/{m_0}}+\Z}t_0^{j} xz^{-j-1}$ on the module $\W$.

 In (\ref{lastimageeqn}), the Virasoro field $\omega_\aff(z)$ is replaced with the twisted 
vertex operator $Y_\W(\omega_\aff,z)$. The latter operator may be written down using (\ref{TBI}) (cf. \cite{KW}):
\begin{eqnarray}
\label{twSug}
Y_\W(\omega_\aff,z) = \frac{1}{2(c+h^\vee)} \bigg(
\sum_i : Y_\W (x_i (-1) \vac, z) Y_\W (x^i (-1) \vac, z) : \nonumber\\
- \sum_i \alpha_i z^{-1}  Y_\W ([x_i, x^i] (-1) \vac, z)
- c \sum_i {{\alpha_i} \choose 2} z^{-2} \Id_\W \bigg) ,
\end{eqnarray}
where $\{x_i\}$, $\{x^i\}$ are dual bases of $\fg$ that are homogeneous relative to the 
grading $\fg = \mathop\bigoplus\limits_{\ol{\bf r}} \fg_{\ol{\bf r}}$ and $\alpha_i$ is a representative 
of the coset $r_0/m_0 +\Z$ for which $x_i \in \fg_{\ol{\bf r}}$.

\subsection{Representations of $\fgtor(\gs_0,\gs_1,\ldots,\gs_N)$}
We are now ready to describe irreducible representations of the twisted toroidal Lie algebra 
$\fgtor(\gs)=\fgtor(\gs_0,\gs_1,\ldots,\gs_N)$.  In this subsection, we describe how 
to construct these representations from the tensor product $\VH \otimes \W \otimes \LglV$ 
of Corollary \ref{1twistmodule}.  We will prove their irreducibility in Section \ref{irreducibility}.

In order to specify the spaces on which $\fgtor(\gs)$ acts, we recall the definition of 
{\em thin covering} of a module \cite{BL}.  Let ${\mathcal L}=\bigoplus_{g\in G}{\mathcal L}_g$ be a Lie algebra graded by a finite abelian group $G$, and let $U$ be a (not necessarily graded) module for $\mathcal{L}$.  A {\em covering} of $U$ is a collection of subspaces $U_g$ $(g\in G)$ satisfying the following axioms
\begin{enumerate}
\item[{\rm (i)}]$\displaystyle{\sum_{g\in G}U_g=U}$
\item[{\rm (ii)}]$\mathcal{L}_gU_h\subseteq U_{g+h}\hbox{\ for all\ }g,h\in G$.
\end{enumerate}
A covering $\{U_g\ |\ g\in G\}$ is a {\em thin covering} if there is no other covering $\{U_g'\ |\ g\in G\}$ of $U$ with $U_g'\subseteq U_g$ for all $g\in G$.

The automorphisms $\gs_1,\ldots,\gs_N$ extend to commuting automorphisms of the twisted affine Lie algebra $\fgaff(\gs_0)$.  This gives a grading of $\fgaff(\gs_0)$ by the finite abelian group 
$\Z^N/\Gamma$. Let $\left\{\W_{\ol{r}}\ |\ \ol{r}\in \Z^N/ \Gamma) \right\}$ be a thin covering of the irreducible bounded 
$\fgaff(\gs_0)$-module $\W$ fixed in Corollary \ref{1twistmodule}.  
The thin coverings of quasifinite modules like $\W$ were classified in \cite{BL}.

\begin{theorem}\label{mainthm} The space
\begin{equation}
\cM = 
\sum\limits_{r\in\Zx^N}  q^r \otimes \FF \otimes  \W_\orr \otimes L_{\gV},
\end{equation} 
is a $\fgtor(\sigma)$-submodule
in $\VH \otimes \W \otimes \LglV$.
\end{theorem}
\proof
We only need to verify that $\cM$ is closed under the action of the twisted fields
$Y_\W (x(-1)\vac, z) K_0(r,z)$ for $x\in \fg_{\ol{\bf r}}$ and $Y_\W(\omega_\aff, z)$. This,
however, follows immediately from the definition of a covering. 
\qed

\section{Irreducibility}\label{irreducibility}

We now state one of the main results of this paper.

\begin{theorem}\label{irrthm} 
Let $\FF$ be the Fock space $\Cx[u_{pj}, v_{pj}\ |\ p=1,\ldots,N,\ j=1,2,\ldots]$. 
Let $\LglV$ be an irreducible highest weight module for the twisted $\glnhat$-Virasoro algebra
with central character given by (\ref{chargammaii}), and let $\W$ be an irreducible bounded module for the twisted affine algebra $\fghat(\gs_0)$ at level $c\neq 0,-h^\vee$.  Let $\left\{ \W_\orr\ |\ \orr \in \Z^N/\Gamma \right\}$ be a thin covering relative to the automorphisms $\gs_1, \ldots, \gs_N$ of $\fgaff(\gs_0)$.  
Then the space
$$\cM = \sum\limits_{r\in\Zx^N}  q^r \otimes \FF \otimes \W_\orr \otimes L_{\gV}$$
is an irreducible module for the twisted toroidal Lie algebra $\fgtor (\gs)$ with the action given
by (\ref{image1eqn}), (\ref{imageiieqn}), (\ref{imageveqn}), (\ref{lastimageeqn}), 
(\ref{twistedaff}), and (\ref{twSug}).
\end{theorem}

 The proof of this theorem will be split into a sequence of lemmas. Consider a nonzero submodule $\cN\subseteq\cM$. We need to show that $\cN = \cM$.

\begin{lemma}\label{uno}
Let $\{ U^A \}$ be the standard monomial basis of the Fock space
$\FF = \Cx[u_{pj}, v_{pj}\ |\ p=1,\ldots,N,\ j=1,2,\ldots]$. 
Then $w = \sum_A U^A \otimes f_A \in \cN$, where
$$f_A \in  \sum\limits_{r\in\Zx^N}  q^r \otimes \W_\orr \otimes L_{\gV}$$
if and only if $1 \otimes f_A \in \cN$ for all $A$.
\end{lemma}

\noindent
\proof
The Lie algebra $\fgtor(\gs)$ contains the components of the fields $k_a(0,z)$,
$\widetilde{d}_a(0,z)$, $a=1,\ldots, N$, which act as multiplication and differentiation operators
on the Fock space $\Cx[u_{pj}, v_{pj}\ |\ p=1,\ldots,N,\ j=1,2,\ldots]$. This Fock space is an irreducible module over the Heisenberg Lie algebra generated by these operators, which proves the claim of this lemma.
\qed

\bigskip

 The key technique for proving irreducibility under the action of some vertex operators is the following observation: 
any subspace stabilized by the moments of the (untwisted) vertex operators $Y(a,z)$ and $Y(b,z)$ is also (setwise) invariant under 
the moments of the vertex operators $Y(a_{(k)}b,z)$ for all $k\in\Z$.  This is an immediate consequence of the 
Borcherds identity \cite{VA} with $k,n \in \Z$:
$$ (a_{(k)} b)_{(n)} = 
\sum\limits_{j\geq 0} (-1)^{k+j+1} { k \choose j }
b_{(n+k-j)} a_{(j)} 
+ \sum\limits_{j\geq 0}  (-1)^j {k \choose j }
a_{(k-j)} b_{(n+j)}. $$ 
The case of twisted modules requires a more delicate analysis.


\begin{lemma}\label{dos}
The space $\cN$ is closed under the action of the vertex operator $Y(\omega_\hyp, z)$.
\end{lemma}
\proof
The Lie algebra fields $k_a(0,z)$ and $\widetilde{d}_a(0,z)$ act as the vertex operators
$Y(u_{a1}, z)$ and $Y(v_{a1},z)$, respectively. Since
$$\omega_\hyp = \sum_{p=1}^N {u_{p1}}_{(-1)} v_{p1},$$
the lemma now follows from the Borcherds identity observation.
\qed

\bigskip

\begin{lemma}\label{tres}
The space $\cN$ is closed under the action of the vertex operators $Y(E_{ab}(-1) \vac, z)$, for all $a,b$.
\end{lemma}
\proof
For $b=1,\ldots,N$ and $r\in\Gamma$, the Lie algebra field $\td_b (r,z)$ is represented by the vertex operator $Y(v_{b1} q^r,z) + \sum_{p=1}^N r_p Y(E_{pb}(-1)q^r,z)$.
Taking $r=0$, this becomes $Y(v_{b1},z)$. Combining this with the fact that
$k_0(r,z)$ is represented by $cY(q^r,z)$, we see that
$\cN$ is invariant under the action of $Y(v_{b1} q^r,z)=:Y(v_{b1},z)Y(q^r,z):$, and hence also under the field $\sum_{p=1}^N r_p Y(E_{pb}(-1)q^r,z)$.
Since $q^{-r}_{(-1)}(E_{pb}(-1)q^r) = E_{pb}(-1) \vac$, we obtain that the space $\cN$ is
invariant under $\sum_{p=1}^N r_p Y(E_{pb}(-1)\vac,z)$.
Finally choosing $r_a = m_a$
and $r_p = 0$ for $p\neq a$, we see that $\cN$ is invariant under the action of
$Y(E_{ab}(-1) \vac, z)$.
\qed

\begin{lemma}\label{quatro}
Let ${\bf r}\in\Z^{N+1}$ and $x\in \fg_{\overline {\bf r}}$. Then $\cN$ is closed under the action of
$q^r Y_\W(x(-1)\vac,z)$.
\end{lemma}
\proof
 The submodule $\cN$ is closed under the action of $Y_\W(x(-1)\vac,z) Y(q^r,z)$ since this 
operator represents the action of the Lie algebra field $x(r,z)$. Let $w\in\cN$. 
We would like to show that the coefficients of all powers of $z$ in
$q^r Y_\W (x(-1)\vac,z) w$ belong to $\cN$. By Lemma \ref{uno}, we may assume that $w$ does not 
involve $u_{pj}, v_{pj}$. In this case,
$$ Y_\W (x(-1)\vac,z) Y(q^r,z) w = q^r Y_\W (x(-1)\vac,z) w + (\hbox{terms involving \ } u_{pj}) .$$
Applying Lemma \ref{uno} again, we conclude that $\cN$ is closed under the action of $q^r Y(x(-1)\vac,z)$.
\qed

\begin{lemma}\label{cinco} The space $\cN$ is closed under the action of the twisted vertex operator $Y_\W(\omega_{\aff},z)$.
\end{lemma}
\proof
The action of $Y_\W(\omega_{\aff},z)$ is given by (\ref{twSug}). The dual bases $\{ x_i \}$, $\{ x^i \}$ may be assumed to 
be homogeneous.  For each $i$, let ${\bf r}^{(i)}\in\Z^{N+1}$ so that $x_i \in \fg_{\ol{{\bf r}}^{(i)}}$ and $x^i \in \fg_{-\ol{{\bf r}}^{(i)}}$.  Then
\begin{align}
: Y_\W( x_i(-1)&\vac, z) Y_\W( x^i(-1)\vac, z) : \nonumber\\
& =\ \ \ \ :\left( q^{r^{(i)}} Y_\W( x_i(-1)\vac, z) \right) \left( q^{-r^{(i)}} Y_\W( x^i(-1)\vac, z) \right) : ,
\end{align}
so Lemma \ref{quatro} implies that $\cN$ is invariant under the operator 
\break
$: Y_\W( x_i(-1)\vac, z) Y_\W( x^i(-1)\vac, z) :$.

Note also that $[x_i, x^i] \in \fg_{\ol{0}}$.  Thus the components of the field 
\break
$\sum_{j \in\Z} t_0^j [x_i, x^i] z^{-j-1}$ belong to $\fghat(\sigma_0)$. Since this field is 
represented by 
\break
$Y_\W ([x_i, x^i](-1) \vac, z)$, we conclude that $\cN$ is invariant under this operator. The last summand in (\ref{twSug}) involves the identity operator, which leaves $\cN$ invariant. This completes the proof of the lemma.
\qed

\begin{lemma}\label{seis} The space $\cN$ is closed under the action of the vertex
operator $Y(\omega_{\gV},z)$.
\end{lemma}
\proof The Lie algebra field $\td_0(0,z)$ is represented by the vertex operator
$$Y(\omega_\hyp,z) + Y(\omega_\gV,z) + Y_\W(\omega_{\aff},z) .$$
Since $\cN$ is closed under $\td_0(0,z)$, Lemmas \ref{dos} and \ref{cinco} imply that $\cN$ is closed under $Y(\omega_{\gV},z)$.
\qed

\bigskip

We are now ready to complete the proof of Theorem \ref{irrthm}. The Fock space $\FF$
is an irreducible module for the Heisenberg subalgebra in $\fgtor(\gs)$ spanned by the components of the fields $\widetilde{d}_a(0,z)$, $k_a(0,z)$, $a=1,\ldots, N$, together with the central element $k_0$.
The space $L_\gV$ is an irreducible module for the twisted $\widehat{\mathfrak{gl}}_N$-Virasoro algebra.


With respect to the commuting automorphisms $\gs_1,\ldots,\gs_N:\ \fgaff(\gs_0)\rightarrow \fgaff(\gs_0)$, we can form the twisted multiloop Lie algebra
$$ L(\fgaff(\gs_0); \gs_1, \ldots, \gs_N) =
\sum_{s\in\Z^N} t^s\otimes\fgaff(\gs_0)_{\ols}.$$
The twisted affine Lie algebra $\fgaff(\gs_0)$ is generated by the subspaces
$t_0^{r_0/m_0}\fg_{\ol{r}}$ for ${\bf r}=(r_0,r)\in\Z^{N+1}$.  The corresponding operators  $q^r Y_\W(x(-1)\vac,z)$, with  $x\in \fg_{\ol {\bf r}}$, thus generate the action of the  Lie algebra $\cL(\fgaff(\gs_0), \gs_1, \ldots, \gs_N)$ on the module
$$ \sum_{r\in\Zx^N} q^{r} \otimes \W_\orr. $$
By \cite[Section 5]{BL}, we see that this space
is a $\Z^{N+1}$-graded-simple module for the twisted multiloop algebra
$\cL(\fgaff(\gs_0), \gs_1, \ldots, \gs_N)$.

To complete the proof of the theorem, we will use the following fact about tensor products of modules.

\begin{lemma}\label{hjk}
Let $A$ and $B$ be associative unital algebras graded by an abelian group $G$.  Suppose that $V$ and $W$ are $G$-graded-simple 
modules for $A$ and $B$, respectively, with $V_\gamma$ finite-dimensional for all $\gamma\in G$. For $\alpha\in G$, denote by $V^{(\alpha)}$ a $G$-graded $A$-module, obtained from $V$ by a shift in grading: $V^{(\alpha)}_\gamma = V_{\gamma+\alpha}$. Assume that for all $\alpha\in G$,
$\alpha \neq 0$, there is no grading-preserving isomorphism of $A$-modules between $V$ and 
$V^{(\alpha)}$. Then $V \otimes W$ is a $G$-graded-simple $A \otimes B$-module. 
\end{lemma}
\proof
 We need to show that $V \otimes W$ can be generated by any homogeneous nonzero element $u \in (V \otimes W)_\gamma$. Let us write 
$$u = \sum\limits_{\alpha\in G} \sum\limits_i v^i_\alpha \otimes w^i_\alpha
\hbox{{\hskip 0.6cm} \rm (finite sum)},$$
where $v^i_\alpha\in V_\alpha, \; w^i_\alpha\in W_{\gamma - \alpha} .$
Without loss of generality, we may assume that the set of vectors $\{ v^i_\alpha \}$ is linearly independent, and the vectors $\{ w^i_\alpha \}$ are nonzero. Let $v \otimes w$
be one of the terms in the sum above. The modules in 
$\left\{ V^{(\alpha)} | \alpha\in G \right\}$ do not have admit grading-preserving isomorphisms between pairs of distinct modules, so we may apply the quasifinite density theorem (\cite[Thm. A.2]{BL}) to conclude that there exists $a \in A_0$ such that 
$a v = v$, while $a v^i_\alpha = 0$ for all other $i$ and $\alpha$.  Thus the $A \otimes B$ submodule generated by $u$ contains $v \otimes w = (a \otimes 1) u$,
where both $v$ and $w$ are homogeneous. Acting on $v \otimes w$ with $A \otimes 1$ and 
$1 \otimes B$, and taking into account that both modules $V$ and $W$ are graded-simple,
we conclude that the submodule generated by $u$ is $V \otimes W$.  
\qed

\bigskip

The space $\FF \ot L_\gV$ is an irreducible module for the universal enveloping algebra
of the direct sum of the infinite-dimensional Heisenberg algebra and the twisted $\glnhat$-Virasoro algebra $\gV$. This module has a natural $\Z$-grading, which we
extend to a $\Z^{N+1}$ grading by setting $\left( \FF \ot L_\gV \right)_{(r_0, r)} = 0$
whenever $r \neq 0$. We immediately see by comparing the characters that there are no 
grading-preserving isomorphisms between this module and the modules obtained from it by shifts in the grading. By the result of Section 5 of \cite{BL}, the space 
$\sum_{r\in\Zx^N} q^{r} \otimes \W_\orr$ is a $\Z^{N+1}$-graded-simple module
for the twisted multiloop algebra $\cL(\fgaff(\gs_0), \gs_1, \ldots, \gs_N)$.
Since the Lie algebra $\fgtor(\gs)$ contains the derivations $d_0, d_1, \ldots, d_N$, the 
space $\cN$ is a $\Z^{N+1}$-graded submodule of $\cM$, and every element of $\cN$ can be reduced to a homogeneous element using $d_0,d_1,\ldots,d_N$.  By Lemmas \ref{dos}--\ref{seis}, $\cN$ is closed under the action of the Heisenberg Lie algebra, of $\gV$, and of $\cL(\fghat(\gs_0), \gs_1, \ldots, \gs_N)$. 
 Applying Lemma \ref{hjk}, we see that $\cN = \cM$, and thus $\cM$ is irreducible.

\section{Irreducible modules for twisted toroidal EALAs}

Irreducible modules for untwisted toroidal extended affine algebras were constructed in \cite{B2}.  The techniques developed in the previous sections can be used to extend this construction and obtain irreducible modules for the twisted toroidal EALAs.

 The twisted toroidal EALA 
$$ \fgE =(\cR\ot\fgdot)\oplus\cK\oplus_\tau\cS $$
is spanned by elements $d_0, d_1, \ldots, d_N$ and
by the moments of the fields
\begin{eqnarray}\label{fielddiv}
k_0(s,z)&=&\sum_{j\in\Z}t_0^{j}t^s k_0z^{-j},\nonumber\\
k_a(s,z)&=&\sum_{j\in\Z}t_0^{j}t^s k_az^{-j-1},\nonumber\\
x(r,z)&=&\sum_{j\in{{r_0}/{m_0}}+\Z}t_0^{j}t^rxz^{-j-1}\hbox{\ for each\ }x\in\fgdot_{\ol{\bf r}},\nonumber\\
\widetilde{d}_{ab}(s,z)&=&\sum_{j\in\Z}\big(s_b t_0^{j}t^s d_a- s_a t_0^{j}t^s d_b \big)z^{-j-1},\nonumber\\
\widehat{d}_a(s,z)&=&\sum_{j\in\Z}\left( j t_0^{j}t^s d_a + s_a t_0^{j}t^s {\widetilde {d}}_0 
+ \frac{s_a}{2cN} (N - 1 + \mu c) t_0^j t^s k_0  \right) z^{-j-2},\nonumber
\end{eqnarray}
where $s \in \Gamma$, ${\bf {r}} \in \Z^{N+1}$, and $a,b = 1, \dots, N$.

 In the representation theory of $\fgE$, the twisted $\glnhat$-Virasoro algebra $\gV$ is replaced 
with its subalgebra, the semidirect product $\sV$ of the Virasoro algebra with the affine algebra $\widehat{\sll}_N$:
$$ \sV = \left( \Cx[t_0, t_0^{-1}] \otimes \sln  \oplus \Cx C_\sln \right) \rtimes 
\left( \Der \Cx[t_0, t_0^{-1}] \oplus \Cx C_\vir \right).$$

 Now we can state the theorem describing irreducible modules for the twisted toroidal EALA
$\fgE (\gs)$.

\begin{theorem} \label{Erep}
Let $\W$ be an irreducible bounded $\fgaff(\gs_0)$-module of level $c \neq 0, - h^\vee$, 
with a thin covering
$\left\{ \W_{\orr}\ |\ \orr \in \Z^N / \Gamma \right\}$ relative to 
automorphisms $\gs_1, \ldots, \gs_N$. Let $L_\sV$ be an irreducible highest weight module for
$\sV$ with a central character $\gamma$:
$$ \gamma (C_\sln) = 1 - \mu c, \; 
\gamma(C_\vir) = 12 (1 - \frac{1}{N}) + 12 \mu c (1 + \frac{1}{N}) 
- 2N - \frac {c \dim (\fg)} {c + h^\vee} .$$
Then the space
$$\sum\limits_{r\in\Zx^N}  q^r \otimes \FF \otimes \W_\orr \otimes L_{\sV} $$
has the structure of an irreducible module for the twisted toroidal extended affine 
Lie algebra $\fgE( \gs_0, \gs_1, \ldots, \gs_N)$ with the action given by
\begin{align}
k_0(s,z)& \mapsto c K_0 (s,z), \\
k_a(s,z)& \mapsto c K_a (s,z), \\
x(r,z)& \mapsto Y_\W (x(-1) \vac, z) K_0 (s,z), \\
\widetilde{d}_{ab}(s,z)& \mapsto :\left( s_b D_a (z) - s_a D_b (z) \right) K_0 (s,z) : \nonumber\\
& + s_b \sum\limits_{{p=1} \atop {p\neq a}} s_p E_{pa}(z) K_0 (s,z)
- s_a \sum\limits_{{p=1} \atop {p\neq b}} s_p E_{pb}(z) K_0 (s,z)\nonumber\\
&+ s_a s_b (E_{aa} - E_{bb})(z) K_0 (s,z), \\
\widehat{d}_a(s,z)& \mapsto s_a : \left( \omega_\hyp(z) + \omega_\sV(z) + 
Y_\W (\omega_\aff,z) \right) K_0 (s,z) : \nonumber\\
+ s_a & \sum\limits_{p,\ell=1}^N s_p \psi_1 (E_{p \ell}) (z) K_\ell (s,z)
+ s_a (\mu c - 1) \sum\limits_{p=1}^N s_p \left( \frac{\partial}{\partial z} K_p(z) \right)
K_0 (s,z) \nonumber\\
 - \bigg( & z^{-1} + \frac{\partial}{\partial z} \bigg)  \left( :D_a(z) K_0(s,z): 
+ \sum\limits_{p=1}^N s_p \psi_1 (E_{pa})(z) K_0 (s,z) \right),  
\end{align}
where $s \in \Gamma$, ${\bf {r}} \in \Z^{N+1}$, $a,b = 1, \dots, N$, and $\psi_1$ is the natural projection (\ref{proj}) 
$\psi_1:\ \gln(\Cx)\rightarrow \sln(\Cx)$.
\end{theorem}

 This theorem is based on its untwisted analogue \cite[Thm. 5.5]{B2}. The proof is completely parallel to the proof of Theorems \ref{mainthm} and \ref{irrthm} and will be omitted.

\

\section{Example: EALAs of Clifford Type}\label{examplessection}

We now apply the general theory that we have developed to construct irreducible 
representations of EALAs coordinatized by Jordan tori of Clifford type.

\subsection{Multiloop realization} 

Fix a positive integer $m$.  Let $2\Z^m\subseteq S\subseteq \Z^m$, where $S$ is a union of some cosets of the subgroup $2\Z^m\subseteq \Z^m$.  We assume that $\Z^m$ is generated by $S$ as a group.  Let $\ol{\mu}$ be the image of $\mu\in\Z^m$ under the map $\Z^m\rightarrow\Z^m/2\Z^m$, and let $r$ be the cardinality of $\ol{S}=S/2\Z^m$.  We identify $\Z^m/2\Z^m$ with the multiplicative group $\{-1,1\}^m=\Z_2^m$.

 A {\em Jordan torus of Clifford type} is a Jordan algebra
$$J=\displaystyle{\bigoplus_{\mu\in S}\Cx s^\mu}$$
with multiplication given by
$$s^\mu s^\eta=\left\{\begin{array}{ll}
s^{\mu+\eta} & \hbox{if}\ \ol{\mu}=\ol{0},\ \ol{\eta}=\ol{0},\ \hbox{or}\ \ol{\mu}=\ol{\eta},\\
0 & \hbox{otherwise}.
\end{array}\right.
$$
Let $L_J=\{L_a\ |\ a\in J\}$ be the set of left multiplication operators $L_a: b\mapsto ab$ on $J$.  The {\em Tits-Kantor-Koecher algebra} associated with $J$ is the Lie algebra
$$\hbox{TKK}(J)=\big(J\ot \mathfrak{sl}_2(\Cx)\big)\oplus [L_J,L_J],$$
where
\begin{align*}
[a\ot  x,b\ot y]&=ab\ot[x,y]\ +\ (x|y)[L_a,L_b],\\
[d,a\ot x]&=da\ot x=-[a\ot x,d],\\
[d,d']&= dd'-d'd
\end{align*}
for all $a,b\in J,\ x,y\in\mathfrak{sl}_2(\Cx)$, and $d,d'\in [L_J,L_J]$.

We now introduce some notation which will be used to realize $\hbox{TKK}(J)$ as a multiloop algebra.
Let $U$ be an $(r+2)$-dimensional vector space with a basis $\{v_i\ |\ i\in \I\}$,
where $\I=\{1,2,3\}\cup\left( \ol{S}\setminus\{\ol{0}\}\right)$.
Define a symmetric bilinear form on $U$ by declaring that this basis is orthonormal.
If $i\in\{1,2,3\}$, let $i_p=1$, and if $i=\ol{\mu}\in\ol{S}\setminus\{\ol{0}\}$,
let $i_p=\ol{\mu}_p\in\{-1,1\}$ for all $p\in\{1,\ldots,m\}$.  For each $i,j,k\in\I$,
define $e_{ij}\in\mathfrak{so}(U)$ by
$$e_{ij}(v_k)=\gd_{jk}v_i-\gd_{ik}v_j.$$
Let $\gs_p$ be an orthogonal transformation on $U$ defined by
$$\gs_p(v_i)=i_pv_i$$
for all $i\in\I$ and $p\in\{1,\ldots,m\}$.  We identify each $\gs_p$ with an automorphism
of $\mathfrak{so}(U)$
where $\gs_p$ acts on $\mathfrak{so}(U)$  by conjugation.
Each of these $\gs_p\in\hbox{Aut}\,(\mathfrak{so}(U))$ has order $2$.

\bigskip

\noindent
{\bf Remark.} In our construction the index set $\I$ is obtained from $\ol{S}$ by triplicating $\ol{0}$
into $\{1,2,3\}$. This is done in order to create a 3-dimensional subalgebra $\mathfrak{so}_3 (\Cx) \cong
\mathfrak{sl}_2 (\Cx)$, fixed under all the involutions $\gs_p$. 
Analogous gradings on $\mathfrak{so}_{2^m} (\Cx)$  were considered in 
\cite{Po,Al,BSZ}.

\begin{theorem}\label{jmulti}
The Tits-Kantor-Koecher algebra $\hbox{TKK}(J)$ is isomorphic to the twisted multiloop algebra $\mathcal{G}=L(\mathfrak{so}_{r+2}(\Cx);\gs_1,\ldots,\gs_m)$ via the following map:
\begin{align*}
\phi:\ \hbox{TKK}(J)&\rightarrow L(\mathfrak{so}_{r+2}(\Cx);\gs_1,\ldots,\gs_m)\\
s^\mu\ot X_1&\mapsto \left\{\begin{array}{ll}
T^\mu\ot e_{32} & \hbox{if}\ \ol{\mu}=\ol{0}\\
T^\mu \ot e_{\ol{\mu} 1} & \hbox{otherwise,}
\end{array}\right.\\
s^\mu\ot X_2 &\mapsto \left\{\begin{array}{ll}
T^\mu \ot e_{13} & \hbox{if}\ \ol{\mu}=\ol{0}\\
T^\mu\ot e_{\ol{\mu} 2} & \hbox{otherwise,}
\end{array}\right.\\
s^\mu\ot X_3&\mapsto \left\{\begin{array}{ll}
T^\mu\ot e_{21} & \hbox{if}\ \ol{\mu}=\ol{0}\\
T^\mu\ot e_{\ol{\mu} 3} & \hbox{otherwise,}
\end{array}\right.\\
[L_{s^\gamma},L_{s^\eta}]&\mapsto T^{\gamma+\eta}\ot e_{\ol{\gamma} \ol{\eta}},
\end{align*}
for all $\mu,\gamma,\eta\in S$ with $\ol{\gamma},\ol{\eta},\ol{\gamma+\eta}\neq\ol{0}$.
Here $\{X_1,X_2,X_3\}$ is a basis of $\mathfrak{sl}_2(\Cx)$ with relations
\begin{eqnarray*}
[X_1,X_2]=X_3,\ [X_2,X_3]=X_1,\ [X_3,X_1]=X_2.
\end{eqnarray*}

\end{theorem}
\proof Observe that $\gs_p(e_{ij})=i_p j_pe_{ij}$ for all $i,j\in\I$ and $p\in\{1,\ldots,m\}$.  This implies that the image of $\phi$ is contained in the twisted multiloop algebra $L(\mathfrak{so}_{r+2}(\Cx);\gs_1,\ldots,\gs_m)$.  The verification that $\phi$ is a homomorphism is tedious but straightforward, and will be omitted.  It is clear that $\phi$ is injective.  To see that it is surjective, we note that $\hbox{TKK}(J)$ and $L(\mathfrak{so}_{r+2}(\Cx);\gs_1,\ldots,\gs_m)$ have natural $\Z^m$-gradings given by
\begin{align*}
\deg(s^\mu \ot X_i)&=\mu,\\
\deg[L_{s^\gamma},L_{s^\eta}]&=\gamma+\eta,\\
\deg(T^\mu\ot e_{ij})&=\mu.
\end{align*}
The map $\phi$ is then homogeneous of degree $0$.  It is now sufficient to verify that the dimensions of the corresponding graded components of $\hbox{TKK}(J)$ and $L(\mathfrak{so}_{r+2}(\Cx);\gs_1,\ldots,\gs_m)$ are the same.

It is easy to see that
\begin{align*}
\dim \hbox{TKK}(J)_\mu&=\dim\hbox{TKK}(J)_\gamma,\\
\dim\mathcal{G}_\mu&=\dim\mathcal{G}_\gamma
\end{align*}
whenever $\ol{\mu}=\ol{\gamma}$.  This allows us to define
\begin{align*}
a'_{\ol{\mu}}&=\dim\hbox{TKK}(J)_\mu\\
a''_{\ol{\mu}}&=\dim\mathcal{G}_\mu
\end{align*}
for all $\mu\in\Z^m$.  Instead of proving that $a'_{\ol{\mu}}=a''_{\ol{\mu}}$ for each $\ol{\mu}\in\Z_2^m$, we will show that
\begin{equation}\label{etoile}
\sum_{\ol{\mu}\in\Z_2^m}a'_{\ol{\mu}}=\sum_{\ol{\mu}\in\Z_2^m}a''_{\ol{\mu}}.
\end{equation}

Since the map $\phi$ is injective and homogeneous of degree $0$, the latter equality will imply that $a'_{\ol{\mu}}=a''_{\ol{\mu}}$, and $\phi$ is thus an isomorphism.  We now verify (\ref{etoile}).

The contribution of the space $J\ot \mathfrak{sl}_2(\Cx)$ in the sum $\sum_{\ol{\mu}\in\Z_2^m}a'_{\ol{\mu}}$ is $3r$, while the space $[L_J,L_J]$ contributes $\binom{r-1}{2}$.  Thus $\sum_{\ol{\mu}\in\Z_2^m}a'_{\ol{\mu}}=3r+{r-1 \choose 2}$.
The right-hand side of (\ref{etoile}) is simply the dimension ${{r+2}\choose{2}}$ of $\mathfrak{so}_{r+2}(\Cx)$.  Since $3r+{{r-1}\choose{2}}={{r+2}\choose{2}}$, we are done.\qed

\bigskip

 We are interested in EALAs associated with the universal central extension of 
$\hbox{TKK}(J)$. The above multiloop realization of $\hbox{TKK}(J)$ yields a description of
such EALAs in the setup of Section 1 as $\fgE(\gs_1,\ldots, \gs_m)$ with $\fg = \so_{r+2} (\Cx)$.

 We now consider the representation theory of these Jordan torus EALAs.
 To conform with the notation in the rest of the paper, we will set $m=N+1$, and we number the variables of the Jordan torus from $0$ to $N$.  Likewise, the automorphisms of $\mathfrak{so}_{r+2}(\Cx)$ under consideration become
$\gs_0,\gs_1,\ldots,\gs_N$ and
the variables in the multiloop algebra are thus changed from $T_1, \ldots, T_m$ to $t_0^{1/2}, t_1,
\ldots, t_N$.

 According to Theorem \ref{mainthm}, the piece of the simple module $\cM$ specific to the Jordan torus EALA is the irreducible highest weight module $\W$ for the twisted affine Lie algebra $\widehat{\mathfrak{so}}_{r+2}(\gs_0)$, and its thin covering 
$\left\{ \W_\orr \right\}$ with respect to the automorphisms
$\{ \gs_1, \ldots, \gs_N \}$.

 The Lie algebra $\widehat{\mathfrak{so}}_{r+2}(\gs_0)$ is isomorphic to the
untwisted or twisted affine Lie algebra, depending on whether or not $\gs_0$ is an inner automorphism of $\mathfrak{so}_{r+2}(\Cx)$.

\begin{lemma}\label{autso}
Let $U$ be a finite-dimensional vector space with an orthonormal basis $\{ v_i\ |\ i\in \I \}$. Let $\gs \in GL(U)$, with $\gs (v_i) = \pm v_i$. Then conjugation by $\gs$ is an inner automorphism of $\mathfrak{so}(U)$ if and only if the matrix of $\gs$ in this basis has an even number of $-1$'s on the diagonal or an even number of $+1$'s.
\end{lemma}
\proof
If $\gs$ has an even number of $-1$'s on the diagonal, then $\gs\in SO(U)$, and it is inner. Suppose that $\gs$ has an even number of $+1$'s on the diagonal. Since
$\gs$ acts on $\mathfrak{so}(U)$ by conjugation, $-\gs$ induces the same automorphism, and $-\gs$ is inner by the previous argument. Finally, suppose that
$\gs$ has an odd number of $+1$'s and an odd number of $-1$'s on the diagonal. Then $\dim U$ is even, and $\mathfrak{so}(U)$ is of type D. Then, multiplying $\gs$
by an appropriate diagonal matrix with an even number of $-1$'s on the diagonal, we can get a diagonal matrix $\tau$ with all $+1$'s except for the entry $-1$ in the last position. Choosing a Cartan subalgebra of $\mathfrak{so}(U)$ and a basis of its root system as in Subsection \ref{FullTKK} below, one can easily see that $\tau$ is the Dynkin diagram automorphism of order 2 of a root system of type $D$. Since $\gs$ differs from $\tau$ by a factor which is an inner automorphism, we conclude that
$\gs$ is not inner.
\qed

\bigskip

 We will now focus our attention on two Clifford type EALAs of nullity 2: a ``baby TKK'', and a ``full lattice TKK''.

\bigskip

\subsection{Baby TKK} 
Let $S$ be a union of 3 cosets of $2\Z^2$ in $\Z^2$,
corresponding to the coset representatives $(0,0),\ (0,1),$ and $(1,0)$.  Then $m=2$ and $r=3$.  The corresponding multiloop algebra given by Theorem \ref{jmulti} is
$L(\mathfrak{so}_5; \gs_0, \gs_1)$ where the matrices
$$\gs_0 =
\left(\begin{array}{lllll}
1&&&&\\
&1&&&\\
&&1&&\\
&&&1&\\
&&&&-1\\
\end{array}\right),
\hbox{\quad\ }
\gs_1 =
\left(\begin{array}{lllll}
1&&&&\\
&1&&&\\
&&1&&\\
&&&-1&\\
&&&&1\\
\end{array}\right)
$$
act on $\mathfrak{so}_5(\Cx)$ by conjugation.  Multiplying by $-I$ does not change their action:

$$\gs_0 =
\left(\begin{array}{lllll}
-1&&&&\\
&-1&&&\\
&&-1&&\\
&&&-1&\\
&&&&1\\
\end{array}\right),
\hbox{\quad\ }
\gs_1 =
\left(\begin{array}{lllll}
-1&&&&\\
&-1&&&\\
&&-1&&\\
&&&1&\\
&&&&-1\\
\end{array}\right).
$$
Consider a Cartan subalgebra in $\mathfrak{so}_5$ with a basis of coroots
$$
h_1 =
\left(\begin{array}{ll|ll|l}
&\quad&\quad&&\;\\
&&&&\\
\hline
&&&2i&\\
&&-2i&&\\
\hline
&&&&\\
\end{array}\right),
\hbox{\quad\ }
h_2 =
\left(\begin{array}{ll|ll|l}
&i&\quad&\quad&\;\\
-i&&&&\\
\hline
&&&-i&\\
&&i&&\\
\hline
&&&&\\
\end{array}\right)
.
$$

 The corresponding generators of the simple root spaces are:
$$e_1 =
\left(\begin{array}{rr|rr|r}
\;&\;&&&\\
&&&&\\
\hline
&&&&-1\\
&&&&i\\
\hline
&&1&-i&\\
\end{array}\right),
\hbox{\quad\ }
f_1 =
\left(\begin{array}{rr|rr|r}
\;&\;&&&\\
&&&&\\
\hline
&&&&1\\
&&&&i\\
\hline
&&-1&-i&\\
\end{array}\right),
$$

$$e_2 =
\frac{1}{2}\left(\begin{array}{rr|rr|r}
&&1&i&\;\\
&&-i&1&\\
\hline
-1&i&&&\\
-i&-1&&&\\
\hline
&&&&\\
\end{array}\right),
\hbox{\quad\ }
f_2 =
\frac{1}{2}\left(\begin{array}{rr|rr|r}
&&-1&i&\;\\
&&-i&-1&\\
\hline
1&i&&&\\
-i&1&&&\\
\hline
&&&&\\
\end{array}\right).
$$

 These generators satisfy the relations $[h_i, e_j] = a_{ij} e_j$, $[h_i, f_j] = -a_{ij} f_j$, $[e_i, f_j] = \delta_{ij}$ with the Cartan matrix
$$A = \left( \begin{array}{rr} 2 & -2 \\ -1 & 2 \\ \end{array}\right).$$
 Note that
$$\gs_0 = \exp\left(\pi i \,\hbox{ad}\,h_2\right).$$

 The Lie algebra $\mathfrak{so}_5$ has an eigenspace decomposition $\mathfrak{so}_5 = \mathfrak{so}_5^{\ol{0}} \oplus \mathfrak{so}_5^{\ol{1}}$ with respect to the action of $\gs_0$:
$$\mathfrak{so}_5^{\ol{j}}=\{x\in\mathfrak{so}_5\ |\ \gs_0x=(-1)^jx\}$$
for $j=0,1$.  Since $\gs_0$ and $\gs_1$ commute, the subspaces  $\mathfrak{so}_5^{\ol{0}}$ and
$\mathfrak{so}_5^{\ol{1}}$ are invariant with respect to $\gs_1$. We view $\gs_1$ as an automorphism
of the twisted loop algebra
$$\cL(\mathfrak{so}_5; \gs_0) = \mathop\sum\limits_{j\in\Z} t_0^{j/2} \mathfrak{so}_5^{\ol{j}},$$
by letting it act by $\gs_1(t_0^{j/2}x)=t_0^{j/2}\gs_1(x)$ for each $x\in\mathfrak{so}_5^{\ol{j}}$ and $j\in\Z$.  We then extend it to an automorphism of the twisted affine algebra
$\widehat{\mathfrak{so}}_5 (\gs_0)$ by $\gs_1(C_\aff) = C_\aff$.

 Since $\gs_0$ is inner (as is every automorphism of $\mathfrak{so}_5$), the twisted loop algebra $\cL(\mathfrak{so}_5; \gs_0)$ is isomorphic to the untwisted loop algebra
$\Cx [t_0, t_0^{-1}]\ot \mathfrak{so}_5 $ \cite[Prop. 8.5]{kac}. This lifts to an isomorphism of affine algebras
$$\theta: \quad \widehat{\mathfrak{so}}_5 \rightarrow \widehat{\mathfrak{so}}_5 (\gs_0),$$
such that
$$\theta(t_0^j\ot e_\alpha) = t_0^{j + \alpha(h_2)/2}\ot e_\alpha ,$$
$$\theta(t_0^j\ot h_\alpha) = t_0^{j}\ot h_\alpha  +
\delta_{j,0} \frac{\alpha(h_2)}{2} C_\aff,$$
$$\theta(C_\aff) = C_\aff,$$
where $e_\alpha$ is in a root space of $\mathfrak{so}_5$ with $\alpha\neq 0$,
and $h_\alpha$ in the Cartan subalgebra is normalized so that $[e_\alpha, e_{-\alpha}] =
(e_\alpha | e_{-\alpha}) h_\alpha$.

 Using the identification $\theta$, we transform $\gs_1$ into an automorphism $\wgs_1 = \theta^{-1} \gs_1 \theta\in\hbox{Aut}\,(\widehat{\mathfrak{so}}_5)$.  Since $\theta$ does not preserve
the natural $\Z$-grading of $\widehat{\mathfrak{so}}_5$, $\wgs_1$ does not leave the components
$t_0^j\ot\mathfrak{so}_5$ invariant. However, it does leave 
invariant the Cartan subalgebra $\mathfrak{h}=(\Cx h_1 \oplus \Cx h_2) \oplus \Cx C_\aff$ of
$\widehat{\mathfrak{so}}_5$.

 Let us describe the group $N$ of automorphisms of an affine Lie algebra $\fgaff$ that leave the Cartan subalgebra
$\mathfrak{h}$ setwise invariant \cite{PK}.  First of all, $N$ has a normal subgroup $H \times \Cx^*$ of automorphisms that fix $\mathfrak{h}$ pointwise.  Here $\Cx^*$ is the set of automorphisms $\{ \tau_a\ |\ a\in\Cx\setminus 0 \}$ that act by
$$ \tau_a (t_0^j\ot x) = a^j t_0^j\ot x, \;\;\; \tau_a (C_\aff) = C_\aff ,$$
and $H= \{ \exp(\ad\,h)\ |\ h \in \h \}$ consists of inner automorphisms. The quotient
$N/(H\times \Cx^*)$ can be presented as follows:
$$ N/(H\times \Cx^*) \cong \langle\pi\rangle \times \left( Aut(\Gamma) \ltimes W \right) ,$$
where $\pi$ is the Chevalley involution, $W$ is the affine Weyl group and $Aut(\Gamma)$ is the group of automorphisms of the affine Dynkin diagram $\Gamma$. The elements of this factor group may be viewed as permutations of the roots of the affine Lie algebra.

 Let $\gs$ be an automorphism of the affine Lie algebra $\fgaff$ leaving invariant its Cartan subalgebra.   Such a Cartan subalgebra may always be found relative to any family of finite order automorphisms $\gs_0,\ldots,\gs_N$ by Appendix \ref{appendix}.  Let $(\W,\rho)$ be an integrable irreducible highest weight module for 
$\fgaff$ with dominant integral highest weight $\lambda$ relative to a fixed base of simple roots.
In order to determine the thin covering of $\W$ with respect to the cyclic group generated by $\gs$, we need to know whether
the modules $(\W,\rho)$ and $(\W,\rho \circ \gs)$ are isomorphic \cite{BL}.
The answer to this question does not change if we replace $\gs$ with $\gs \circ \mu$, where
$\mu$ is an automorphism of $\fgaff$ for which $(\W,\rho)$ and $(\W,\rho \circ \mu)$
are isomorphic.

\begin{proposition}\label{propX} Let $(V, \rho)$ be a representation of a Lie algebra $L$. Suppose $x\in L$, $\ad\,x$ is locally nilpotent on $L$, and $\rho(x)$ is locally nilpotent on $V$. Then the representations $(V,\rho)$
and $(V, \rho \circ \exp (\ad\,x))$ are isomorphic.
\end{proposition}
\proof 
 The isomorphism from $(V,\rho)$ to $(V, \rho \circ \exp (\ad\,x))$ is given by
$\exp (\rho(x))$. This claim is equivalent to the identity 
$$ \exp(\rho(x)) \rho(y) = \rho (\exp((\ad\,x)y)) \exp (\rho(x)), \; y\in L,$$
which is well-known.
\qed

\begin{corollary} \label{inner}
Let $(\W,\rho)$ be an integrable module for an affine Lie algebra $\fgaff$, and let
$\mu$ be an inner automorphism of $\fgaff$. Then the modules $(\W,\rho)$ and $(\W,\rho \circ \mu)$ are isomorphic.
\end{corollary}
\proof
 This follows from the previous proposition and the fact that the Kac-Moody group (of inner automorphisms) is generated by the exponentials of the real root elements, which are locally nilpotent on integrable modules
\cite{PK}.
\qed

\begin{lemma} \label{scale}
 Let $(\W,\rho)$ be an integrable irreducible highest weight module for $\fgaff$.
Let $a\in\Cx^*$.
Then the modules $(\W,\rho)$ and  $(\W,\rho \circ \tau_a)$ are isomorphic.
\end{lemma}
\proof
 It is easy to see that the module $\W$ admits a compatible action of the group $\Cx^*$
(see \cite[Section 4]{PK}), which we will denote by
$$ T_a : \W \rightarrow \W, \;\;\; a\in\Cx^* .$$
This can be done by requiring that $\Cx^*$ fixes the highest weight vector and satisfies the compatibility condition
$$T_a \rho(y) v = \rho(\tau_a (y)) T_a v .$$
The condition implies that $T_a$ is a module isomorphism between $(\W,\rho)$ and  $(\W,\rho \circ \tau_a)$.
\qed

\bigskip

 Let $\bgs$ be the image of $\gs$ in the factor-group 
$N/W(H \times \Cx^*) \cong \langle\pi\rangle \times Aut(\Gamma)$. 
We will also identify $\bgs$ with an automorphism of $\fgaff$ by viewing
$Aut(\Gamma)$ as a subgroup of $Aut(\fgaff)$.

 Since $WH$ consists of inner automorphisms of $\fgaff$, Corollary \ref{inner} and Lemma
\ref{scale} imply the following lemma:

\begin{lemma}\label{hatbar}
Let $(\W,\rho)$ be an integrable irreducible highest weight module for $\fgaff$ and $\gs$ be an
automorphism of $\fgaff$ leaving invariant its Cartan subalgebra.
Then the $\fgaff$-modules $(\W,\rho \circ \gs)$  and $(\W,\rho \circ \bgs)$ are isomorphic.\qed
\end{lemma}

We now return to the setting of our particular example and calculate 
$\bgs_1 = \overline{\theta^{-1} \gs_1 \theta}$ 
as an  automorphism of the root system of $\widehat{\mathfrak{so}}_5$. The Dynkin diagram $\Gamma$ of
$\widehat{\mathfrak{so}}_5$ is

\medskip

\centerline{
\renewcommand{\tabcolsep}{1mm}
\begin{tabular}{p{1mm}p{4mm}p{1mm}p{4mm}p{1mm}}
\renewcommand{\tabcolsep}{1mm}
$\alpha_0$ && $\alpha_1$ && $\alpha_2$ \\
$\bullet$ & $\Rightarrow$ &$\bullet$ & $\Leftarrow$ & $\bullet$ \\
\end{tabular}
}

\medskip

\noindent
and $Aut(\Gamma) \cong \Z_2$.

 We need to compute the action induced by $\wgs_1=\theta^{-1} \gs_1 \theta$ on the simple roots $\alpha_0, \alpha_1, \alpha_2$ of $\widehat{\mathfrak{so}}_5$. Since $\wgs_1$ leaves invariant the null root spaces, the induced automorphism of the root system fixes the null root $\delta$. Taking into
account that $\alpha_0 = \delta - 2\alpha_1 - \alpha_2$, we see that it is enough to
find the action on $\alpha_1$ and  $\alpha_2$:
$$ \wgs_1 (e_1) = \theta^{-1} \gs_1 \theta (e_1) =
\theta^{-1} t_0^{\alpha_1 (h_2)/2}\ot \gs_1 (e_1) $$
$$ = \theta^{-1} t_0^{-1/2}\ot \gs_1
\left(\begin{array}{rr|rr|r}
&&&&\;\\
&&&&\\
\hline
&&&&-1\\
&&&&i\\
\hline
&&1&-i&\\
\end{array}\right)= \theta^{-1} t_0^{-1/2}\ot
\left(\begin{array}{rr|rr|r}
&&&&\;\\
&&&&\\
\hline
&&&&-1\\
&&&&-i\\
\hline
&&1&i&\\
\end{array}\right)$$
$$=\theta^{-1} (- t_0^{-1/2}  \ot f_1) = -t_0^{-1/2 + \alpha_1(h_2)/2}\ot f_1
 =  -t_0^{-1}\ot f_1  .$$
Thus we get
$\wgs_1(\alpha_1) = - \delta - \alpha_1$.
Similarly,
$$ \wgs_1 (2e_2) = \theta^{-1} \gs_1 \theta (2e_2) = \theta^{-1} t_0^{\alpha_2 (h_2)/2}\ot
 \gs_1 (2e_2) $$
$$ = \theta^{-1} t_0\ot \gs_1
\left(\begin{array}{rr|rr|r}
&&1&i&\;\\
&&-i&1&\\
\hline
-1&i&&&\\
-i&-1&&&\\
\hline
&&&&\\
\end{array}\right)=\theta^{-1} t_0\ot
\left(\begin{array}{rr|rr|r}
&&1&-i&\;\\
&&-i&-1&\\
\hline
-1&i&&&\\
i&1&&&\\
\hline
&&&&\\
\end{array}\right)$$
$$= \theta^{-1} t_0\ot [[e_1, e_2], e_1] = t_0^{1 - (2\alpha_1(h_2)+\alpha_2(h_2))/2} \ot [[e_1, e_2], e_1]
 =  t_0\ot [[e_1, e_2], e_1].$$
We get
$\wgs_1(\alpha_2) = \delta + 2\alpha_1 + \alpha_2$, and hence
$\wgs_1(\alpha_0) = 2\delta - \alpha_2$.
Let $\gamma$ be the diagram automorphism that interchanges $\alpha_0$ with $\alpha_2$ and fixes $\alpha_1$. Let $r_0, r_1, r_2$ be the simple reflections generating the affine Weyl group of
$\widehat{\mathfrak{so}}_5$. It is then straightforward to verify that as an automorphism of the root
system,
$$\wgs_1 = r_1 r_2 r_0 r_1 \gamma .$$
We conclude that $\bgs_1 = \gamma \in Aut(\Gamma)$.

\begin{proposition}
Let $(\W,\rho)$ be an irreducible highest weight module for $\widehat{\mathfrak{so}}_5$ of dominant
integral highest weight $\lambda$.
The $\widehat{\mathfrak{so}}_5$-modules $(\W,\rho)$ and $(\W,\rho \circ \wgs_1)$
are isomorphic if and only if the diagram automorphism $\gamma$ fixes the highest weight $\lambda$.
\end{proposition}
\proof
 Taking into account the fact that $\bgs_1 = \gamma$ and applying Lemma \ref{hatbar}, we
conclude that the modules $(\W,\rho)$ and $(\W,\rho \circ \wgs_1)$ are isomorphic if and only if
the modules $(\W,\rho)$ and $(\W,\rho \circ \gamma)$ are isomorphic. However, $(\W,\rho \circ \gamma)$
is the highest weight module with the highest weight $\gamma(\lambda)$. Since two irreducible highest
weight modules are isomorphic precisely when their highest weights are equal, we obtain the claim of the
proposition.
\qed

\begin{corollary}
If the highest weight $\lambda$ is not fixed by the diagram automorphism $\gamma$, then the thin 
covering of $\W$ with respect to the cyclic group $\langle\wgs_1\rangle \cong \Z_2$ is
$\{ \W, \W \}$.
\end{corollary}
\proof
This follows from Theorem 4.4 of \cite{BL} and the proposition above.
\qed

\bigskip

When the highest weight $\lambda$ {\em is} fixed by the diagram automorphism $\gamma$, there is a module isomorphism $\phi_{\gamma}:\ (\W,\rho)\rightarrow (\W,\rho\circ\gamma)$.  Concretely, we may define the action of $\phi_\gamma$
on the Verma module of highest weight $\lambda$ by postulating that $\phi_\gamma$ fixes the highest 
weight vector and $\phi_\gamma(\rho(x)v)=\rho(\gamma(x))\phi_\gamma(v)$, for all $x\in\widehat{\mathfrak{so}}_5$ and $v\in\W$. It is also clear that $\phi_\gamma$ will leave invariant the maximal submodule of this 
Verma module. This gives an action of $\gamma$ as the operator $\phi_\gamma$ on $\W$. As a result, we obtain an action
of the semi-direct product $\langle\gamma\rangle \ltimes G$ on $\W$, where $G$ is the Kac-Moody group of $\widehat{\mathfrak{so}}_5$. 
Since $\wgs_1 \in \langle\gamma\rangle \ltimes G$, we realize $\wgs_1$ as an order 2 operator on $\W$. It is
easy to see that the action of $\wgs_1$ on $\W$ is locally finite, and hence $\W$ has 
decomposition 
\begin{equation}
\label{smallcovering}
\W = \W_{\overline 0} \oplus  \W_{\overline 1}, 
\end{equation}
where $\W_{\overline 0}$, $\W_{\overline 1}$ are the $\pm 1$ eigenspaces of $\wgs_1$. 
In this case, $\left\{ \W_{\overline 0}, \W_{\overline 1} \right\}$ is a thin covering of $\W$ relative to
$\wgs_1$.  We have now proved the following theorem:

\begin{theorem}\label{baby}
Let $(\W, \rho)$ be an irreducible highest weight module for $\widehat{\mathfrak{so}}_5$ of integral dominant 
highest weight $\lambda$. 
 View $\W$ as a module for the twisted affine 
algebra $\widehat{\mathfrak{so}}_5 (\gs_0)$ with the action $\rho \circ \theta^{-1}$. 
Let $L_{\sV}$ be an irreducible highest weight module for $\sV$ as in Theorem \ref{Erep}. 

(i) If  $\gamma(\lambda) \neq \lambda$, where $\gamma$ is the Dynkin diagram automorphism of 
$\widehat{\mathfrak{so}}_5$
then the space
$$\VH \otimes L_{\sV} \otimes \W$$
has the structure of an irreducible module for the ``Baby TKK'' EALA $\fgE (\gs_0, \gs_1)$
with the action described in Theorem \ref{Erep}.  

(ii) If $\gamma(\lambda) = \lambda$, then the space
$$\sum\limits_{r\in\Zx} q^{r} \otimes \FF \otimes L_{\sV} \otimes  \W_\orr $$
is an irreducible module for $\fgE (\gs_0, \gs_1)$, where $\left\{ \W_{\overline 0}, \W_{\overline 1} \right\}$ are as in (\ref{smallcovering}).\qed
\end{theorem}

\bigskip

\subsection{Full lattice TKK of nullity 2} 
\label{FullTKK}
Let $S$ be the set of all 4 cosets of $2\Z^2$ in $\Z^2$.
The corresponding multiloop algebra given by Theorem \ref{jmulti} is
$L(\mathfrak{so}_6; \gs_0, \gs_1)$ where
$$\gs_0 =
\left(\begin{array}{llllll}
1&&&&&\\
&1&&&&\\
&&1&&&\\
&&&1&&\\
&&&&-1&\\
&&&&&-1\\
\end{array}\right),
\hbox{\quad\ }
\gs_1 =
\left(\begin{array}{llllll}
1&&&&&\\
&1&&&&\\
&&1&&&\\
&&&-1&&\\
&&&&1&\\
&&&&&-1\\
\end{array}\right).
$$

 We consider a Cartan subalgebra in $\mathfrak{so}_6$ with a basis of coroots
$$
h_1 =
\left(\begin{array}{ll|ll|ll}
&\quad&\quad&&&\\
&&&&&\\
\hline
&&&i&&\\
&&-i&&&\\
\hline
&&&&&i\\
&&&&-i&\\
\end{array}\right),
\hbox{\quad\ }
h_2 =
\left(\begin{array}{ll|ll|ll}
&i&\quad&\quad&\quad&\\
-i&&&&&\\
\hline
&&&-i&&\\
&&i&&&\\
\hline
&&&&&\\
&&&&&\\
\end{array}\right),$$
$$h_3 =
\left(\begin{array}{ll|ll|ll}
&\quad&\quad&&&\\
&&&&&\\
\hline
&&&i&&\\
&&-i&&&\\
\hline
&&&&&-i\\
&&&&i&\\
\end{array}\right),
\hbox{\quad\ }
\label{hiii}$$

 The corresponding generators of the simple root spaces are:
$$e_1 =
\frac{1}{2}
\left(\begin{array}{rr|rr|rr}
&\quad&&&&\\
&&&&&\\
\hline
&&&&1&-i\\
&&&&-i&-1\\
\hline
&&-1&i&&\\
&&i&1&&\\
\end{array}\right),
\hbox{\quad\ }
f_1 =
\frac{1}{2}
\left(\begin{array}{rr|rr|rr}
&\quad&&&&\\
&&&&&\\
\hline
&&&&-1&-i\\
&&&&-i&1\\
\hline
&&1&i&&\\
&&i&-1&&\\
\end{array}\right),
$$

$$e_2 =
\frac{1}{2}\left(\begin{array}{rr|rr|rr}
&&1&i&\quad&\\
&&-i&1&&\\
\hline
-1&i&&&&\\
-i&-1&&&&\\
\hline
&&&&&\\
&&&&&\\
\end{array}\right),
\hbox{\quad\ }
f_2 =
\frac{1}{2}\left(\begin{array}{rr|rr|rr}
&&-1&i&\quad&\\
&&-i&-1&&\\
\hline
1&i&&&&\\
-i&1&&&&\\
\hline
&&&&&\\
&&&&&\\
\end{array}\right),
$$

$$e_3 =
\frac{1}{2}
\left(\begin{array}{rr|rr|rr}
&\quad&&&&\\
&&&&&\\
\hline
&&&&1&i\\
&&&&-i&1\\
\hline
&&-1&i&&\\
&&-i&-1&&\\
\end{array}\right),
\hbox{\quad\ }
f_3 =
\frac{1}{2}
\left(\begin{array}{rr|rr|rr}
&\quad&&&&\\
&&&&&\\
\hline
&&&&-1&i\\
&&&&-i&-1\\
\hline
&&1&i&&\\
&&-i&1&&\\
\end{array}\right).
$$

 The root system of $\mathfrak{so}_6$ is of type $D_3 = A_3$ and
 the Lie brackets $[h_i, e_j] = a_{ij} e_j$, $[h_i, f_j] = -a_{ij} f_j$ are given by the Cartan matrix
$$A = \left( \begin{array}{rrr} 2 & -1 & 0 \\ -1 & 2 & -1 \\ 0 & -1 & 2 \\ \end{array}\right).$$
 Note that
$$\gs_0 = \exp\left(\hbox{ad}\,\frac{\pi}{2} i (h_1 - h_3)\right).$$

 Since $\gs_0$ is an inner automorphism of $\mathfrak{so}_5$, the twisted loop
algebra $\cL(\mathfrak{so}_6; \gs_0)$ is isomorphic to the untwisted loop algebra
$\Cx [t_0, t_0^{-1}]\ot\mathfrak{so}_6  $ \cite[Prop. 8.5]{kac}.
This lifts to an isomorphism $\theta$ of affine Lie algebras:
$$\theta: \quad \widehat{\mathfrak{so}}_6 \rightarrow \widehat{\mathfrak{so}}_6 (\gs_0),$$
where
$$\theta(t_0^j\ot e_\alpha) = t_0^{j + \alpha(h_1-h_3)/4}\ot e_\alpha ,$$
$$\theta(t_0^j\ot h_\alpha) =  t_0^{j}\ot h_\alpha +
\delta_{j,0} \frac{\alpha(h_1 - h_3)}{4} C_\aff,$$
$$\theta(C_\aff) = C_\aff .$$

 Using the identification $\theta$, we transform $\gs_1$ into an automorphism $\wgs_1 = \theta^{-1} \gs_1 \theta$ of $\widehat{\mathfrak{so}}_6$.  Let $(\W,\rho)$ be an integrable irreducible highest weight module for $\widehat{\mathfrak{so}}_6$ with dominant integral highest weight $\lambda$.  We now determine when
the modules $(\W,\rho)$ and $(\W,\rho \circ \wgs_1)$ are isomorphic.

Let $\bgs_1$ be the image of $\wgs_1$ in $\langle\pi\rangle \times Aut(\Gamma)$ under the projection
$N \rightarrow N/W(H \times \Cx^*)$.  The Dynkin diagram $\Gamma$ of
$\widehat{\mathfrak{so}}_6$ is
%


{\center
{\vskip -0.2cm}
{\hskip 00pt}
$\alpha_0$
{\hskip 20pt}
$\alpha_3$
{\vskip -3pt}
{\hskip 03pt}
$\bullet$
{\hskip -8pt}
${\rule [2pt] {32pt} {1pt}}$
{\hskip -8pt}
$\bullet$
{\vskip -2pt}
{\hskip 03pt}
${\rule {1pt} {32pt}}$
{\hskip 28pt}
${\rule {1pt} {32pt}}$
{\vskip -12pt}
{\hskip 03pt}
$\bullet$
{\hskip -8pt}
${\rule [2pt] {32pt} {1pt}}$
{\hskip -8pt}
$\bullet$
{\vskip -5pt}
{\hskip 00pt}
$\alpha_1$
{\hskip 20pt}
$\alpha_2$
{\vskip 0.5cm}
}

\noindent
so $Aut(\Gamma)$ is the dihedral group of order 8.

As in the case of the Baby TKK algebra described above, we view $\wgs_1$ as an automorphism of the affine root system and calculate the action of $\wgs_1$ on the simple roots 
$\alpha_0, \alpha_1, \alpha_2, \alpha_3$ of $\widehat{\mathfrak{so}}_6$.
After making the analogous calculations, we obtain
$$\wgs_1 = r_3 r_0 r_2 r_1 \eta,$$
where $\eta$ is the diagram automorphism switching $\alpha_0$ with $\alpha_2$, and $\alpha_1$
with $\alpha_3$.

\begin{proposition}
Let $(\W,\rho)$ be an irreducible highest weight module for $\widehat{\mathfrak{so}}_6$ with a dominant
integral weight $\lambda$.
If $\eta(\lambda) \neq \lambda$, then the thin covering of $\W$ with respect to the cyclic group $\langle\wgs_1\rangle$ is
$\{ \W, \W \}$.\qed
\end{proposition}

When $\eta(\lambda)=\lambda$, the module $\W$ and its twist by $\widehat{\gs_1}$ are isomorphic, as in the case of the Baby TKK algebra.  The isomorphism defines a $\Cx$-linear action of $\wgs_1$ as an order 2 operator on $\W$, which then composes into eigenspaces $\W_{\ol{0}}$ and $\W_{\ol{1}}$ relative to the action.  We thus obtain the analogue of Theorem \ref{baby}:
\begin{theorem}
Let $(\W, \rho)$ be an irreducible highest weight representation of $\widehat{\mathfrak{so}}_6$ of integral dominant highest weight $\lambda$. 

(i)  Suppose that $\eta(\lambda) \neq \lambda$, where $\eta$ is the $180^\circ$ rotation automorphism of the Dynkin diagram of $\widehat{\mathfrak{so}}_6$. View $\W$ as a module for the twisted affine algebra $\widehat{\mathfrak{so}}_6 (\gs_0)$ with the action $\rho \circ \theta^{-1}$. Let $L_{\sV}$ be an irreducible highest weight module for $\sV$ as in Theorem \ref{Erep}. Then the space
$$\VH \otimes L_{\sV} \otimes \W$$
has the structure of an irreducible module for the nullity 2 full lattice TKK EALA $\fgE (\gs_0, \gs_1)$
with the action described in Theorem \ref{Erep}.

(ii)  If $\eta(\lambda)=\lambda$, then the space 
$$\sum\limits_{r\in\Zx} q^{r} \otimes \FF \otimes L_{\sV} \otimes  \W_\orr $$
is an irreducible module for $\fgE (\gs_0, \gs_1)$ under the action described in Theorem \ref{Erep}.
\qed  
\end{theorem}

\appendix\label{appendix}

%
%
%

\section{Appendix}
Let $\gs_0,\ldots,\gs_N$ be commuting finite-order automorphisms of a finite-dimensional simple Lie algebra $\fg$.  In this appendix, we show that there exists a  Cartan subalgebra $\fh\subseteq\fg$ which is (setwise) invariant under these automorphisms.  

\begin{lemma}\label{invariantform}
Let $\gs$ be an automorphism of a finite-dimensional simple Lie algebra $L$ with Killing form $(\cdot\,|\,\cdot)$.  Then $(\gs x|\gs y)=(x|y)$ for all $x,y\in L$.
\end{lemma}
\proof Fix elements $x,y\in L$, set $f=\ad x\,\ad y$, and set $g=\ad(\gs x)\,\ad(\gs y)$.  Let $\{v_i\}$ be a basis of $L$.  For each $z\in L$, write
$$f(z)=\sum_if_i(z)v_i$$
with $f_i(z)\in F$.  Then
$$(x|y)= \hbox{tr}(f)=\sum_if_i(v_i).$$
Since $\gs$ is an automorphism, $\{\gs v_i\}$ is also a basis of $L$.  With respect to this basis, we see that
\begin{align*}
g(\gs v_i)&=[\gs x,[\gs y,\gs v_i]]\\
&=\gs[x,[y,v_i]]\\
&=\gs f(v_i)\\
&=\sum_j f_j(v_i)\gs v_j.
\end{align*}
Therefore,
$$(\gs x|\gs y)=\hbox{tr}(g)=\sum_if_i(v_i),$$
and $(x|y)=(\gs x|\gs y)$.\qed

\bigskip

The following lemma appears as \cite[Lm 8.1]{kac} in the context of simple Lie algebras.  The fact that $L^\gs$ is reductive had also appeared previously in \cite[\S1, no. 5]{bourbaki} and \cite[Chap III]{jacobson}.

\begin{lemma}\label{kaclemma}
Let $\gs$ be a finite-order automorphism of a finite-dimensional reductive Lie algebra $L$.  Let $H$ be a Cartan subalgebra of the fixed point subalgebra $L^\gs=\{x\in L\ |\ \gs x=x\}$.  Then $L^\gs$ is reductive, and the centralizer $\Cx_L(H)$ of $H$ in $L$ is a Cartan subalgebra of $L$.
\end{lemma}
\proof The same arguments given in \cite[Lm 8.1]{kac} hold for the reductive case as well.  The only exception is the justification for the $\gs$-invariance of the Killing form, which we have already verified in Lemma \ref{invariantform}.\qed

\bigskip

We can now prove the existence of a $\gs_1,\ldots,\gs_N$-invariant Cartan subalgebra:

\begin{theorem}
Let $\gs_1,\ldots ,\gs_N$ be commuting finite-order automorphisms of a finite-dimensional reductive Lie algebra $L$.  Then $L$ has a (setwise) $\gs_1,\ldots,\gs_N$-invariant Cartan subalgebra $\fh$.
\end{theorem}
\proof  We induct on the dimension of $L$.  If the dimension of $L$ is $0$ or $1$, then the theorem holds since $\fh=L$.  
We can also assume that $\gs_1:\ L\rightarrow L$ is not the identity map.

By Lemma \ref{kaclemma}, the ($\gs_1,\ldots,\gs_N$-invariant) fixed point subalgebra
$$L^{\gs_1}=\{x\in L\ |\ \gs_1 x=x\}$$
is reductive, and by the induction hypothesis, it has a $\gs_1,\ldots,\gs_N$-invariant Cartan subalgebra $H$.  By Lemma \ref{kaclemma}, the centralizer $\fh=\Cx_L(H)$ in $L$ is a Cartan subalgebra of $L$.  It is also $\gs_1,\ldots,\gs_N$-invariant, since
$$[\gs_ix,h]=\gs_i[x,\gs_i^{-1}(h)]\in \gs_i[x,H]=0$$
for all $x\in\fh$, $h\in H$, and $i=1,\ldots,N$.\qed

\end{document}